\documentclass[12pt]{article}

\usepackage{amssymb}
\usepackage{graphicx}
\usepackage[utf8]{inputenc}
\usepackage{amsmath}
\usepackage{amsthm}
\usepackage{amssymb}
\usepackage[colorlinks=true]{hyperref}
\usepackage{subfig}
\usepackage[left=2.5cm,right=2.5cm]{geometry}
\usepackage[english]{babel}

\usepackage[toc,page]{appendix} 
\usepackage{tabularx} 

\usepackage{tikz}
\usetikzlibrary{shapes,arrows}

\usepackage{authblk}

\theoremstyle{plain}

\begin{document}

\title{Adaptive multi-stage integrators for optimal energy conservation in molecular simulations}

\author[1]{Mario Fern\'andez-Pend\'as \thanks{Contact address: \texttt{mfernandez@bcamath.org}}}
\author[1,2]{Elena Akhmatskaya}
\author[3]{J. M. Sanz-Serna}
\affil[1]{BCAM - Basque Center for Applied Mathematics, Alameda de Mazarredo 14, E-48009 Bilbao, Spain}
\affil[2]{IKERBASQUE, Basque Foundation for Science, E-48013 Bilbao, Spain}
\affil[3]{Departamento de Matem\'aticas, Universidad Carlos III de Madrid, Avenida de la Universidad 30, E-28911 Legan\'es (Madrid), Spain}

\maketitle

\begin{abstract}
  We introduce  a new Adaptive Integration Approach (AIA) to be used in a wide range of  molecular simulations. Given a simulation problem and a step size, the method automatically chooses the optimal scheme out of an available family of numerical integrators. Although we focus on two-stage splitting integrators, the idea may be used with more general families.
  In each instance, the system-specific integrating scheme identified by our approach is optimal in the sense that it provides the best conservation of energy for harmonic forces.
  The AIA method has been implemented in the BCAM-modified GROMACS software package.
  Numerical tests in molecular dynamics  and hybrid Monte Carlo simulations of constrained and unconstrained physical systems show that the method successfully realises the fail-safe strategy.  In all experiments, and for each of the criteria employed,  the AIA is at least as good as, and often significantly outperforms the standard Verlet scheme, as well as fixed parameter, optimized two-stage integrators.
  In particular, for the systems tested where harmonic forces play an important role, the sampling efficiency found in simulations using the AIA is up to 5 times better than the one achieved with other tested schemes.
\end{abstract}

\section{Introduction}\label{intro}
  We introduce  a new Adaptive Integration Approach (AIA) to be used in a wide range of  molecular simulations. Given a molecular simulation problem and a step size $\Delta t$, the method automatically chooses the optimal scheme out of an available family of numerical integrators. Although we focus on two-stage splitting integrators \cite{BlanesCasasSanzSerna}, the idea may be used with more general families.
  The system-specific integrating scheme identified by our approach is optimal in the sense that it provides the best conservation of energy for harmonic forces. For hybrid Monte Carlo (HMC) methods \cite{Duane,Horowitz,Kennedy,curso}, the chosen scheme may be expected to provide the biggest possible acceptance rate in the Metropolis accept-reject test.

  The efficiency, and even the feasibility, of molecular dynamics (MD) simulations depend crucially on the choice of numerical integrator. The Verlet algorithm is currently the method of choice; its algorithmic simplicity and optimal stability properties make it very difficult to beat, as discussed in \cite{BlanesCasasSanzSerna}. Splitting integrators may however offer the possibility of improving on Verlet, at least in some circunstances. Those integrators evaluate the forces more than once per step and, due to their simple kick-drift structure, may be implemented easily by modifying existing implementations of the  Verlet scheme. Here we study two-stage integrators. There is a one-parameter family of them \cite{BlanesCasasSanzSerna}, and the parameter value that results in a method with smallest error constant was first identified by McLachlan \cite{MinimunError}. While McLachlan's scheme is the best choice in any given problem if the step length $\Delta t$ is very small, it turns out that its stability interval is not long. This entails that in molecular simulations, where  small time steps are prohibitively expensive, McLachlan's method is likely not to be a good choice. One has then to sacrifice the size of the error constant to ensure that the integrator is able to operate satisfactorily with larger step sizes. Recommended in \cite{BlanesCasasSanzSerna} is a parameter value that achieves a balance between good conservation of energy for reasonable values of $\Delta t$ and  accuracy for small $\Delta t$. That parameter value does not vary with the problem being considered or with the value of $\Delta t$ attempted by the user. On the contrary, in the AIA suggested here,
  the parameter value is automatically adjusted for each problem and each choice of $\Delta t$. On stability grounds, for any given problem, there is a maximum possible value of $\Delta t$; beyond this maximum all integrators in the family are unstable. When the step size chosen by the user is near the maximum value, AIA picks up an integrator that is  (equivalent to) the standard Verlet scheme. As $\Delta t$ decreases, AIA  changes the integrator to ensure optimal conservation of energy;  for $\Delta t$ close to 0, AIA chooses McLachlan's scheme.

  The ideas behind the method are presented in \ref{Formulation}. We also explain how to extend the algorithm to cases with holonomic constraints. As described in \ref{Implementation}, we have implemented the AIA in the BCAM modified GROMACS software package \cite{GROMACSPaper95,GROMACSPaper08}; this modification \cite{Escribano,FernandezPendasNPT} was developed to achieve better accuracy and sampling performance by means of the incorporation of hybrid Monte Carlo methods and multi-stage numerical integrators.
  Section \ref{sec:Testing} presents the  problems used to
  test the performance of the novel adaptive scheme in molecular dynamics  and HMC simulations of constrained and unconstrained physical systems. Section \ref{Results} is devoted to numerical results.
  The performance of the AIA method is compared with the standard velocity Verlet algorithm and the two-stage integrators with the fixed parameter values suggested in \cite{BlanesCasasSanzSerna} and  \cite{Predescu}.  In all experiments and for each of the criteria employed, the performance of AIA is at least as good as, and often significantly better than,  the performances of the Verlet scheme and the fixed parameter two-stage integrators. The final section presents our conclusions.

\section{Adaptive Integration Approach}\label{Formulation}
  The section provides the formulation of the algorithm suggested in this paper.

  \subsection{The one-parameter family of two-stage integrators}
    We consider Hamiltonians  $H$ that can be written as a sum $H = A + B$ of two partial Hamiltonian functions
    \begin{equation*}
      A = \frac{1}{2} p^T M^{-1} p, \qquad
      B = V(q),
    \end{equation*}
    that  respectively correspond to the kinetic and potential energies; $q$ denotes the positions, $p$ the momenta and $M$ is the mass matrix.
    The equations of motion associated with $H$ are
    \begin{equation}\label{eq:motion}
      \begin{split}
	\frac{\mathrm{d}}{\mathrm{d} t} q &= \nabla_p A(q,\, p) = M^{-1} p, \\
	\frac{\mathrm{d}}{\mathrm{d} t} p &= -\nabla_q B(q,\, p) = -\nabla_q V(q).
      \end{split}
    \end{equation}

    For the partial Hamiltonians $A$ and $B$ the equations of motion may of course be integrated in closed form. In fact, for $A$ the solution is a {\em drift} in position
    \begin{equation*}
      (q(t),\, p(t)) = \varphi_t^A (q(0),\, p(0)),\quad q(t) = q(0) + t M^{-1} p(0),\quad p(t) = p(0),
    \end{equation*}
    and for $B$ the solution is a momentum {\em kick}
    \begin{equation*}
      (q(t),\, p(t)) = \varphi_t^B (q(0),\, p(0)),\quad q(t) = q(0),\quad p(t) = p(0) - t \nabla_q V(q(0)).
    \end{equation*}
    Here $\varphi_t^A$ and $\varphi_t^B$ denote the exact solution flows of the partial systems, i.e., the maps that associate with each initial condition $(q(0),\, p(0))$ the exact solution value $(q(t),\, p(t))$.

    The integration schemes under study belong to the family of two-stage splitting methods of the form \cite{BlanesCasasSanzSerna}
    \begin{equation}\label{eq:IntegratorDef2}
      \psi_{\Delta t} = \varphi_{b {\Delta t}}^B \circ \varphi_{{\Delta t}/2}^A \circ \varphi_{(1 -2 b) {\Delta t}}^B \circ \varphi_{{\Delta t}/2}^A \circ \varphi_{b {\Delta t}}^B.
    \end{equation}
    Here $b$ is a parameter, $0<b<1/2$, that identifies the particular integrator being considered and $\psi_{\Delta t}$ denotes the mapping that advances the numerical solution over one step of length ${\Delta t}$.\footnote{It would be possible to consider \lq position\rq\ integrators obtained by swapping the symbols $A$ and $B$ in (\ref{eq:IntegratorDef2}); however the present study just uses the \lq velocity\rq\ form (\ref{eq:IntegratorDef2}). } Note that $\psi_{\Delta t}$ is symplectic \cite{ssc}  as the composition of symplectic mappings and it is time-reversible as a consequence of the palindromic structure of (\ref{eq:IntegratorDef2}).
    The transformation $\varPsi = \varPsi_{{\Delta t},I}$ that advances the numerical solution over $I$ steps is given by the composition
    \begin{equation*}
      \varPsi = \varPsi_{{\Delta t},I} = \overbrace{\psi_{\Delta t} \circ \psi_{\Delta t} \circ \cdots \circ \psi_{\Delta t}}^{I \text{\ times}}.
    \end{equation*}

    Even though $\varphi^B$ appears three times in (\ref{eq:IntegratorDef2}), the methods essentially require {\em two} evaluations of the force $-\nabla V_q$ per step: the  evaluation implicit in the leftmost $\varphi_{b {\Delta t}}^B$ in (\ref{eq:IntegratorDef2}) at the current step is reused in the rightmost $\varphi_{b {\Delta t}}^B$ at the next step. Hence the terminology {\em two-stage} integrator. A fair comparison, in terms of computational cost, between an integration consisting of $I$ steps of length $\Delta t$ with a method of the form (\ref{eq:IntegratorDef2}) and an integration with the standard Verlet integrator uses Verlet  with $2I$ steps of length $\Delta t/2$ (which, in view of Verlet being second order accurate, provides errors that are roughly $1/4$ of those given by Verlet with $I$ steps of length $\Delta t$).

    It is useful in what follows to note that \eqref{eq:IntegratorDef2} may be rewritten as
    \begin{equation}\label{eq:rewrite}
      \psi_{\Delta t} = \Big(\varphi_{b {\Delta t}}^B \circ \varphi_{{\Delta t}/2}^A \circ \varphi_{(1/2 - b) {\Delta t}}^B \Big)
      \circ
      \Big(\varphi_{(1/2 - b) {\Delta t}}^B \circ\varphi_{{\Delta t}/2}^A \circ \varphi_{b {\Delta t}}^B\Big).
    \end{equation}
    The map $\varphi_{(1/2 - b) {\Delta t}}^B \circ\varphi_{{\Delta t}/2}^A \circ \varphi_{b {\Delta t}}^B$ advances the solution over a first half step of length $\Delta t/2$ and is followed by the map $\varphi_{b {\Delta t}}^B \circ \varphi_{{\Delta t}/2}^A \circ \varphi_{(1/2 - b) {\Delta t}}^B$ that effects a second half step, also of length
    $\Delta t/2$. In the particular case  $b=1/4$ both of these maps correspond to a step of length $\Delta t/2$ of the velocity Verlet algorithm:
    \begin{equation*}
      \psi_{{\Delta t}} = \Big(\varphi_{{\Delta t}/4}^B \circ \varphi_{{\Delta t}/2}^A \circ \varphi_{{\Delta t}/4}^B\Big) \circ \Big(\varphi_{{\Delta t}/4}^B \circ \varphi_{{\Delta t}/2}^A \circ \varphi_{{\Delta t}/4}^B\Big)= \varPsi_{{\Delta t}/2}^{VV} \circ\varPsi_{{\Delta t}/2}^{VV}.
    \end{equation*}
    For other values of $b$ the half step maps in \eqref{eq:rewrite} do not coincide with the map of the velocity Verlet integrator, because the durations  $b\Delta t$ and
    $(1/2-b)\Delta t$ are different from one another.
    However, regardless of the choice of $b$, the half step maps have the same kick/drift/kick structure of  velocity Verlet.  This makes it easy
    to implement \eqref{eq:rewrite} by modifying  software that implements the Verlet scheme: it is mainly a matter of adjusting the durations of kicks and drifts (see section \ref{Implementation}).

  \subsection{Nonadaptive choices of the parameter $b$}
    Let us now discuss how best to choose the value of $b$. Regardless of the value of  $b$ the method is second order accurate, i.e.\ the size of the error over one step may be bounded by $C{\Delta t}^3+ \mathcal{O}(\Delta t^5)$, where $C>0$ varies with $b$. McLachlan \cite{MinimunError} was the first to point out that the minimum error constant $C$ is achieved when $b \approx 0.1932$; this is then the optimal value in the limit $\Delta t\rightarrow 0$ of very small step lengths. In molecular dynamics, simulations with values of ${\Delta t}$ that are too small (relatively to the time scales present in the problems) are often unfeasible due to their cost; one may aim to operate with large values of ${\Delta t}$, provided that they are not so large  that the integrations become  unstable. Unfortunately the minimum error constant method possesses a short stability interval $(0,2.55)$ and therefore may not be the best choice when $\Delta t$ is large.
    The stability  of (\ref{eq:IntegratorDef2}) is maximized \cite{BlanesCasasSanzSerna} when $b = 1/4$ with a stability interval
    $(0,4)$. We recall that, for this value of the parameter, integrations with (\ref{eq:IntegratorDef2}) are really Verlet integrations with time step $\Delta t/2$, hence, for $b=1/4$, the stability interval of (\ref{eq:IntegratorDef2}) is twice as long as the stability interval $(0,2)$ of Verlet \cite{Schlick,Skeel,SkeelZhangSchlick}). In fact, it is well known, see e.g.~\cite{BlanesCasasSanzSerna}, that among all explicit integrators that use $k$ force evaluations per step, the longest possible stability interval is obtained by concatenating $k$ Verlet substeps each of length $\Delta t/k$.

    From the discussion above we conclude that the most useful range of values of $b$ is $0.1932 \leq b\leq 0.2500$.  As $b$ increases in this range, both the error constant and the stability interval increase, {\em thus trading accuracy for small $\Delta t$ by the possibility of running stably with larger values of $\Delta t$}.

    The paper \cite{BlanesCasasSanzSerna} recommends the intermediate value $b  \approx 0.2113$. Let us review the ideas leading to this choice, as they will be used in the derivation of our adaptive approach.
    Considered in \cite{BlanesCasasSanzSerna} is the use of algorithms of the form (\ref{eq:IntegratorDef2}) for hybrid Monte Carlo and related simulations. There and in other  situations, the aim is to minimize the energy error
    \begin{equation*}
      \Delta(q,\, p,\, {\Delta t}) = H(\varPsi_{{\Delta t},I}(q,\, p)) - H(q,\, p).
    \end{equation*}
    The analysis in \cite{BlanesCasasSanzSerna} focuses on the model problem where the potential energy is quadratic (harmonic forces), which corresponds to Gaussian probability distributions.
    With the help of a change of variables, the study of the model problem may be reduced to that of the standard harmonic oscillator in nondimensional  variables (standard univariate Gaussian):
    \begin{equation}  \label{eq:stahar}
      (d/dt) q = p, \qquad
      (d/dt) p = - q.
    \end{equation}
    Assume then that the problem (\ref{eq:stahar}) is integrated  by means of (\ref{eq:IntegratorDef2}) and, for reasons that will be apparent immediately, denote by $h$ the step size.
    The expectation or average  $\mathbb{E}(\Delta)$ of the energy error over all possible initial conditions is shown in \cite{BlanesCasasSanzSerna} to possess the bound
    \begin{equation*}
      0 \leq \mathbb{E}(\Delta) \leq \rho(h, b),
    \end{equation*}
    where
    \begin{equation*}
      \rho(h,\, b) = \frac{h^4 (2 b^2 (1/2 - b) h^2 + 4 b^2 -6 b + 1)^2}{8 (2 - b h^2) (2 - (1/2 - b) h^2) (1 - b (1/2 - b) h^2)}.
    \end{equation*}
    It is understood that $\rho=\infty$ for combinations of $b$ and $h$ leading to a denominator $\leq 0$; these combinations correspond to unstable integrations. Thus choices of $b$ and $h$ that lead to a small value of $\rho$ will result in small energy errors for (\ref{eq:stahar}). The study of the function $\rho$ is {\em more discriminating} than the study of the stability interval of the integrators: it is possible for two integrators to share a common stability interval and yet have very different values of $\rho$ for a given value of $h$ that is stable for both of them.

    Let us now move from the scalar oscillator (\ref{eq:stahar})  to multidimensional linear oscillatory problems integrated with step length $\Delta t$ and denote by $\omega_j$, $j = 1,2, \dots$ the corresponding angular frequencies (the periods are $T_j = 2\pi/\omega_j$). By superposing the different modes of the solution, one sees that if the (nondimensional) quantities
    $h_j = \omega_j\Delta t = 2\pi \Delta t/T_j$ are such that, as $j$ varies, all the values $\rho(h_j,b)$ are small, then the energy errors will also be small.
    In \cite{BlanesCasasSanzSerna}, the authors aimed to identify {\em one} value of $b$ that would result in small values of $\rho(h,\, b)$ over a meaningful  range of values of $h$. More precisely, the recommended $b = 0.2113$ was found by minimizing the function of $b$ given by
    \begin{equation}\label{eq:Norm2}
      \max_{0 <  h <  2} \rho(h, b).
    \end{equation}
    The range $0< h < 2$  was chosen because, for the test problems  considered, the standard Verlet method was found to perform well for  $0 < \omega_j\Delta t=2\pi\Delta t/T_j< 1$ (which is half the maximum allowed by the Verlet linear stability interval $(0,2)$); since, as we emphasized already,  (\ref{eq:IntegratorDef2}) uses two force evaluations per step and Verlet only one,  for (\ref{eq:IntegratorDef2}) to be an improvement on standard Verlet it must be demanded that it works well for twice as long values of $\Delta t$, i.e. for $0 < \omega_j\Delta t=2\pi\Delta t/T_j<2$.

    Numerical tests in \cite{BlanesCasasSanzSerna} show the merit of the choice $b = 0.2113$. However the fact remains
    that, if, for a given problem and $\Delta t$, the maximum of $ \omega_j\Delta t=2\pi\Delta t/T_j$ as $j$ varies is significantly smaller than 2, i.e., the chosen $\Delta t$ is relatively small, then  a smaller value of $b$ would provide a better integrator. On the other hand, if that maximum is significantly larger than 2, then it would be advisable to increase $b$.

    A different approach is taken in the present study. Rather than choosing a single value of $b$ that is later applied in all simulations, we suggest an algorithm that, once the system to be integrated has been specified and the user has chosen a value of $\Delta t$, identifies the \lq best\rq\ $b$.

  \subsection{Adapting  the integrator to the problem}
    Although the physical systems that one wishes to simulate in practice are very complex, it is helpful to consider the case  where the forces are two-body interactions. Note that the most stringent stability restrictions on $\Delta t$ are likely to stem from stiff two-body forces, in particular from pairs of bonded atoms. For relatively small energy values, those stiff forces may be assumed to be harmonic.

    For two particles attracting each other harmonically, the period of the oscillations is
    \begin{equation}\label{eq:harmonicperiod}
      T = 2\pi\sqrt{\frac{\mu}{k}},\qquad \mu = \frac{m_1m_2}{m_1+m_2},
    \end{equation}
    where $m_1$, $m_2$ are the masses of the particles, $\mu$ the reduced mass and $k$ the force constant.

    The stability of the integration is of course determined by the highest frequency
    $\tilde{\omega}$ or, equivalently, the smallest period $\tilde T$ present in the system. For the standard Verlet integrator, the linear stability restriction is,
    as noted above,
    \begin{equation}\label{eq:LinearStability2}
      \Delta t < \frac{2}{\tilde \omega} =  \frac{\tilde{T}}{\pi}.
    \end{equation}
    Due to nonlinear effects, including nonlinear resonances, and to other difficulties \cite{SanzSerna91,MandziukSchlick,SchlickMandziukSkeelSrinivas,Skeel,Schlick}, this requirement may be too weak to ensure stability in practice.
    Some authors suggest that the stability restriction for the Verlet integrator
    \begin{equation}\label{eq:ResonanceStability}
      \Delta t < \frac{\sqrt{2}}{\tilde{\omega}} = \frac{\tilde{T}}{\sqrt{2}\pi};
    \end{equation}
    is more realistic in applications than (\ref{eq:LinearStability2}) \cite{Skeel,Schlick,Mazur}.
    Note that moving from (\ref{eq:LinearStability2}) to (\ref{eq:ResonanceStability}) may be seen as the result of multiplying the smallest period by a safety factor $1/\sqrt{2}$ (equivalently multiplying the frequency by $\sqrt{2}$).

    In our adaptive method, if $\Delta t$ is the step size attempted by the user, we exploit the stability restriction in (\ref{eq:ResonanceStability}) to form, similarly to the preceding section, the nondimensional quantity
    \begin{equation}\label{eq:hbar}
      \bar h = \sqrt{2}\tilde{\omega}\Delta t= \sqrt{2}\frac{2\pi}{\tilde{T}}\Delta t
    \end{equation}
    and determine $b$ so as to minimize (cf. (\ref{eq:Norm2}))
    \begin{equation}\label{eq:Normbar}
      \max_{0 <  h <  \bar{h}} \rho(h, b).
    \end{equation}
    Here the function $\rho$ that bounds the energy error is minimized in the {\em shortest interval} $(0,\bar{h})$ that contains all the values
    $ \sqrt{2}\omega_j\Delta t$, where $\omega_j$ are the frequencies in the problem being integrated. Let us illustrate how this works. If the user attempts a value of $\Delta t$ slightly smaller than $ \sqrt{2} \tilde T /\pi$,
    then $\bar h$  will be just below  4 and the minimization of (\ref{eq:Normbar}) will lead to $b$ close to $0.25$. For this value of $b$, $I$ steps of length $\Delta t$ are, as discussed above, equivalent to $2I$ steps of length $\Delta t =   \tilde T /\sqrt{2}\pi$ of the Velocity Verlet algorithm; in other words the adaptive algorithm will run the optimally stable Verlet with the maximum $\Delta t$ allowed by (\ref{eq:ResonanceStability}). As the value  of $\Delta t$ attempted by the user decreases from $ \sqrt{2} \tilde T /\pi$ towards $0$, the value of $b$ will decrease from $0.2500$ to McLachlan's  $0.1932$, thus improving the error constant. The length of the stability interval will shrink
    as $b$ is decreased, but this will  cause no problem because by construction all values $\omega_j\Delta t$ will fall in the stability interval
    (in fact, for safety, even the larger $\sqrt{2}\omega_j\Delta t$ will lie on the stability interval). Finally if $ \Delta t \geq \sqrt{2} \tilde T /\pi$, the quantity (\ref{eq:Normbar}) will be $\infty$ for all values of $b$; this indicates that $\Delta t$ is too large for the problem at hand.

  \subsection{Algorithm}
    Given a physical system and a value of $\Delta t$, the AIA algorithm determines the value of the parameter $b$ to be used in \eqref{eq:IntegratorDef2} as follows:
    \begin{enumerate}
      \item Use equation (\ref{eq:harmonicperiod}) to find the periods or frequencies of all two-body interactions in the system. Determine the minimum period $\tilde T$ and compute the nondimensional quantity  $\bar h$ in \eqref{eq:hbar}.
      \item Check whether $\bar h < 4$. If not, there is no value of $b$ for which the scheme \eqref{eq:IntegratorDef2} is stable for the attempted step size $\Delta t$ and the integration is aborted. In other case go to the next step.
      \item Find the optimal value of the parameter $b$ by minimizing (\ref{eq:Normbar}) with the help of an optimization routine.
    \end{enumerate}

  \subsection{Extension to constrained dynamics}\label{Ap:Constraints}
    Holonomic constraints $g(q) = 0$ allow the use of bigger time steps in physical systems that contain high frequency modes. By freezing those modes, it is possible to bypass the demanding restriction they would otherwise impose on the time step. SHAKE  \cite{SHAKE} and  RATTLE \cite{RATTLE}
    are widely used algorithms in this connection. We now show how, by following the idea behind RATTLE,  two-stage integrators of the family (\ref{eq:IntegratorDef2}) may be applied to problems with constraints. In this way the adaptive integration approach may be extended to the constrained case.

    The constrained equations of motion corresponding to (\ref{eq:motion}) are
    \begin{align*}
      \frac{\mathrm{d}}{\mathrm{d} t} q &= M^{-1} p, \\
      \frac{\mathrm{d}}{\mathrm{d} t} p &= -\nabla_q V(q) + g'(q)^T \lambda, \\
      g(q) &= 0,
    \end{align*}
    where  $\lambda$ is the vector of Lagrange multipliers and $g'(q)^T \lambda$ represents the forces exerted by the constrains. The holonomic constraint implies, by differentiation with respect to time, a constraint on the velocities $(d/dt)q = M^{-1}p$:
    \begin{equation*}
      g^\prime(q)M^{-1}p = 0.
    \end{equation*}
    As in \eqref{eq:rewrite}, we divide one step into two half steps. The equations for the first are
    \begin{equation*}
      \begin{split}
	  p_{n + b} &= p_n - b h \nabla_q V(q_n) + b h g'(q_n)^T \lambda_n,\\
	q_{n + 1/2} &= q_n + \frac{h}{2} M^{-1}p_{n+b},	
      \end{split}
    \end{equation*}
    where the Lagrange multiplier $\lambda_n$ is chosen to ensure
    \begin{equation*}
      g(q_{n+1/2})=0,
    \end{equation*}
    and
    \begin{equation*}
      p_{n + 1/2} =p_{n + b} - (\frac{1}{2} -  b) h \nabla_q V(q_{n + 1/2}) + (\frac{1}{2} -  b) h g'(q_{n + 1/2})^T \lambda_{n + 1/2}^{(v)},
    \end{equation*}
    where the velocity Lagrange multiplier $\lambda_{n + 1/2}^{(v)}$ is chosen so that
    \begin{equation*}
      g'(q_{n+1/2})M^{-1}p_{n+1/2}=0.
    \end{equation*}
    The equations for the second half step $(q_{n+1/2}, p_{n+1/2})\rightarrow (q_{n+1}, p_{n+1})$ are similar.
    The proof of the symplecticness of RATTLE given in \cite{LeimkuhlerSkeel} may be easily adapted to prove that each half step $(q_{n}, p_{n})\rightarrow (q_{n+1/2}, p_{n+1/2})$, $(q_{n+1/2}, p_{n+1/2})\rightarrow (q_{n+1}, p_{n+1})$
    is symplectic. Therefore the whole step $(q_{n}, p_{n})\rightarrow (q_{n+1}, p_{n+1})$ is also symplectic.

    It is clear that  $2I$ steps of  length $\Delta t/2$ of the Verlet integrator supplemented with the constraining technique envisaged here are as expensive as $I$ steps of lenght $\Delta t$  of the extension of two-stage schemes to constrained dynamics we have just described.

    Hybrid Monte Carlo methods can be easily used in constrained dynamics. Only one consideration has to be made: right after the Metropolis test, when the momenta $p_{\small{\mbox{new}}}$ are resampled, the constraint $g'(q) M^{-1}p_{\small{\mbox{new}}} = 0$ has to be fulfilled. Further details can be found in \cite{HMCConstrained}.

\section{Implementation}\label{Implementation}
  \subsection{MultiHMC-GROMACS}
    AIA has been implemented in the MultiHMC-GROMACS software code, an in-house modified version of GROMACS. GROMACS \cite{GROMACSPaper95,GROMACSPaper08} is a popular software package  for molecular dynamics simulations of systems consisting of hundreds to millions of particles, such as proteins, lipids or nucleic acids.  GROMACS supports  state-of-the-art molecular dynamics algorithms and offers extremely fast calculation of non-bonded atomic interactions, which usually are the dominant part of molecular dynamics simulations. It is mainly written in C, highly parallelized, optimized and distributed under the GPL license.

    Currently MultiHMC-GROMACS is based on GROMACS 4.5.4 \cite{GROMACSPaper13}  though its migration to later versions of GROMACS, to take an advantage of CUDA-based GPU acceleration on GPUs \cite{GROMACSPaper15}, is underway.

    MultiHMC-GROMACS has been developed to achieve better accuracy and sampling performance in GROMACS through the use of hybrid Monte Carlo methods and multi-stage numerical integrators. The new algorithms introduced in GROMACS via MultiHMC-GROMACS are the following:
    \begin{itemize}
      \item Hybrid Monte Carlo (HMC) \cite{Duane}, Generalized Hybrid Monte Carlo (GHMC) \cite{Horowitz,Kennedy}, Generalized Shadow Hybrid Monte Carlo (GSHMC) \cite{AkhmatskayaReich,Patent}.

      The implementation of GSHMC in GROMACS has been discussed in detail in  \cite{Escribano,FernandezPendasNPT}. HMC and GHMC are implemented as special cases of GSHMC  \cite{AkhmatskayaReich}.
      Additional parameters have been introduced in the standard GROMACS input file \textit{.mdp} to initiate the new functionalities. These parameters as well as their optional values are presented in the following fragment of \textit{.mdp}:
      {\small
      \begin{verbatim}
	; Hybrid Monte Carlo =
	method                  = HMC;       HMC / GHMC / GSHMC / No
	parameter_phi           = 0.2;       0<phi<pi/2
	nr_MD_steps             = 1000;      any integer
	canonical_temperature   = 310;       any rational
	momentum_flip           = yes;       yes / no
      \end{verbatim}
      }
      \item Multi-stage (two-, three- and four-stage) integrators for unconstrained dynamics.

      The routine \textit{do\_md()} in \textit{md.c}, which performs the integration of the equations of motion, is modified in such a manner that velocity Verlet  steps are concatenated in different ways to form various multi-stage integrators exploiting the Trotter nature of the original implementation \cite{GROMACSPaper13}.
      The parameters needed to construct the desired integrator among all members of the multi-stage families are defined only once, at the beginning of the simulation, and passed to the \textit{update\_coords()} routine in \textit{update.c}, where the actual integration is performed.
      The implementation of multi-stage integrators in MultiHMC-GROMACS is general enough to allow the use of all  members of the families introduced in \cite{BlanesCasasSanzSerna} (we denote them two-s, three-s, four-s, etc.) and in Predescu et al. \cite{Predescu} (two-s-HOH, etc.) as well as the minimum error integrator from \cite{MinimunError} (two-s-minE).  The integrators resulting from the Adaptive Integration Approach described above (two-s-AIA) belong to the family \eqref{eq:IntegratorDef2} of \cite{BlanesCasasSanzSerna} and thus are naturally included in the list.

      No extra variables are needed in the \textit{.mdp} input file. All currently available values of the variable \lq\lq integrator" in \textit{.mdp} are shown below:
      {\small
	\begin{verbatim}
	  integrator = md / md-vv / two-s / two-s-AIA / two-s-HOH /
	               two-s-minE / three-s / four-s
      \end{verbatim}
      }
      \item Two-stage integrators for constrained dynamics.

      The SHAKE algorithm is implemented in the released version of GROMACS using the original approach in \cite{SHAKE}, combined with the Lagrange multipliers procedure of \cite{Lippert2007}  for improving the accuracy in the calculation of velocities of constrained particles \cite{GROMACSPaper08}.

      The implementation of the RATTLE step in GROMACS is done following the algorithm in \cite{RATTLE}.
      The modifications explained in \autoref{Ap:Constraints} for the two-stage integrators for constrained dynamics are combined with the implementation of the released version of GROMACS.
      Any further developments in terms of performance, parallelization or formulation have  not been considered so far.
      \item The v-rescale thermostat \cite{Bussi2007} for two-stage integrators.

      The routine \textit{update\_coupling()} is adapted for two-stage schemes in such a way that the rescaling of the velocities is accurately performed.
      \item AIA  has been introduced in the preprocessing module of GROMACS.
      Its implementation will be described in detail presently.
    \end{itemize}

  \subsection{Implementing AIA}
    \begin{figure}[ht!]
      \centering
      \tikzstyle{decision} = [diamond, draw, fill=white!20,text width=5.5em, text badly centered, node distance=2cm, inner sep=0pt]
      \tikzstyle{block} = [rectangle, draw, fill=white!20,text width=14em, rounded corners, minimum height=4em]
      \tikzstyle{block2} = [rectangle, draw, fill=white!20,text width=20em, node distance=5cm, rounded corners, minimum height=4em]
      \tikzstyle{block22} = [rectangle, draw, fill=white!20,text width=16em, node distance=5cm, rounded corners, minimum height=4em]
      \tikzstyle{block3} = [rectangle, draw, fill=white!20,text width=5em, text badly centered, rounded corners, minimum height=2em]
      \tikzstyle{block4} = [rectangle, draw, fill=white!20,text width=25em, rounded corners, node distance=2cm, minimum height=3em]
      \tikzstyle{line} = [draw, -latex']
      \tikzstyle{cloud} = [draw, ellipse,fill=white!20, node distance=8cm,minimum height=2em]

      \begin{tikzpicture}[node distance = 3cm, auto]
	\node [block] (init) {\textbf{Input}\\
	                      * Modified \textit{.mdp} file\\
	                      * Standard GROMACS input};
 	\node [cloud, right of=init] (system) {\textbf{\textit{Adaptivity?}}};
 	\node [block2, below of=init] (identify) {\textbf{Adaptivity}\\
						  1. Do 1. as in \textbf{No adaptivity} case\\
 	                                          2. Set up the limits for the parameter $b$:
	
	                                          \ \ \ \ $b \in [b_1,\, b_2] \equiv {\cal B}$ \cite{BlanesCasasSanzSerna}\\
 	                                          3. Calculate the fastest period $\tilde{T}$ \eqref{eq:harmonicperiod} and
	
	                                          \ \ \ \ the dimensionless time step $\bar{h}$ \eqref{eq:hbar}.
	
	                                          \ \ \ \ Take $h \in (0,\, \bar{h}) \equiv {\cal H}$\\
 	                                          4. For each $b \in {\cal B}$ calculate $\displaystyle\max_{h \in {\cal H}}{\rho(h,\, b)}$ \eqref{eq:Normbar}\\
 	                                          5. Find optimal $b$ as $\displaystyle\arg \min_{b \in {\cal B}} \max_{h \in {\cal H}}{\rho(h,\, b)}$\\
	
 	                                          6. Pass value of \lq integrator' and optimal $b$\\
						  \ \ \ \ to \textit{.tpr}};
 	\node [block22, below of=system] (identify2) {\textbf{No adaptivity}\\
 	                                             1. For all pairs of particles:\\
 	                                             \ \ \ \ 1.a. Calculate  period $T$ in \eqref{eq:harmonicperiod}\\
 	                                             \ \ \ \ 1.b. If $5 \Delta t \leq T$, STOP\\
 	                                             \ \ \ \ 1.c. If $10 \Delta t \leq T$, WARNING\\
 	                                             2. Pass value of \lq integrator' to \textit{.tpr}};
 	\node [block3, yshift=-9.5cm, xshift=3.5cm] (evaluate) {\textit{.tpr} file};
  	\node [block4, below of=evaluate] (decide) {* Define the integrator in the Trotter factorization form \\
	                                            * Run MD};
 	\path [line] (system) -- node {yes} (identify);
 	\path [line] (system) -- node {no} (identify2);
 	\path [line] (identify) -- (evaluate);
	\path [line] (identify2) -- node {\scriptsize \textbf{Preprocessed input}} (evaluate);
 	\path [line] (evaluate) -- node {\scriptsize \textbf{Runner (\textit{mdrun})}} (decide);
 	\path [line] (init) -- node[align=center] {\scriptsize $\textbf{Preprocessor}$\\ \scriptsize$\textbf{(\textit{grompp})}$} (system);
      \end{tikzpicture}
      \caption{Flowchart of the Adaptive Integration Approach (AIA) as implemented in GROMACS.}\label{fig:flowChartAdaptive}
    \end{figure}
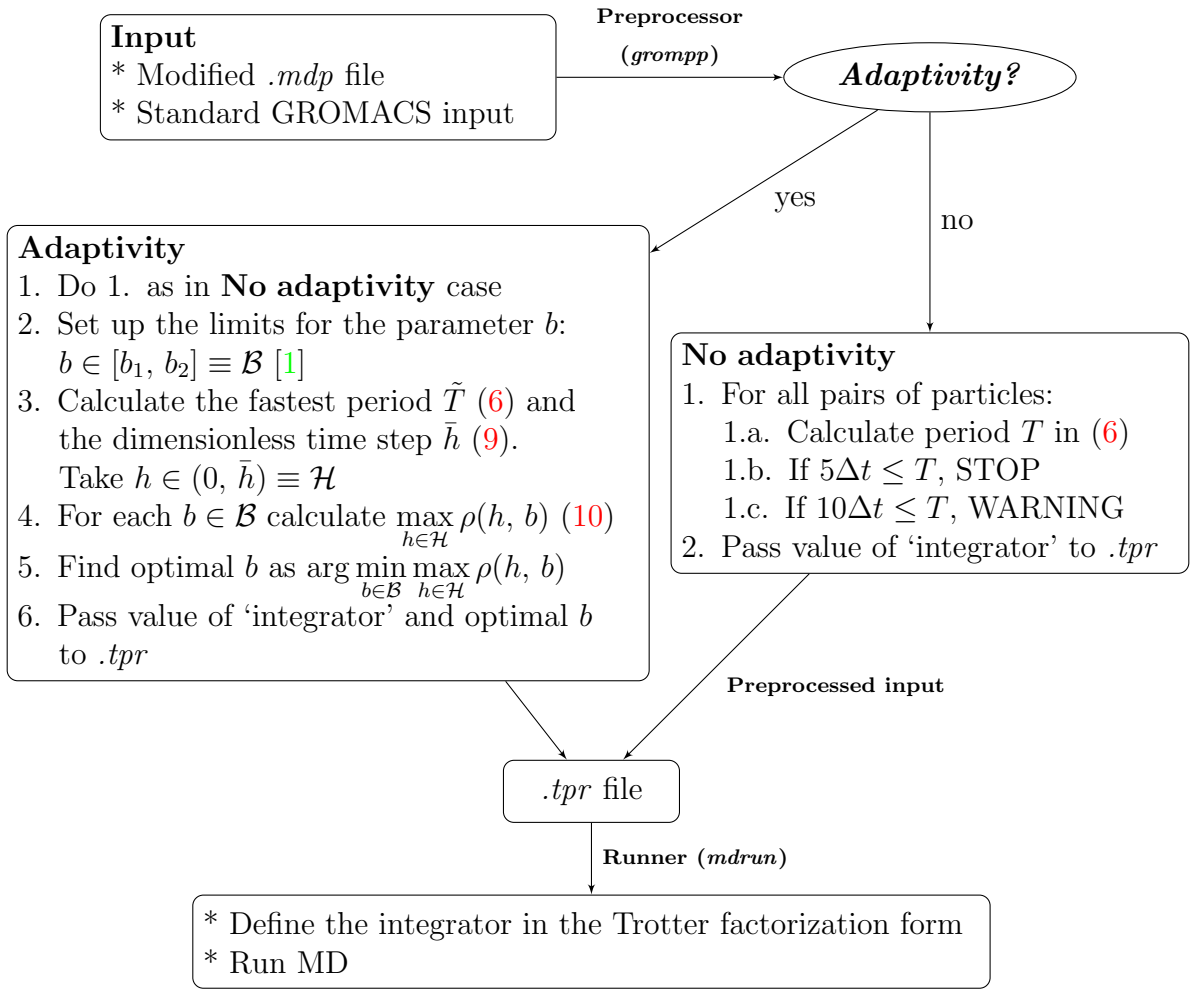

    AIA has been implemented in the GROMACS preprocessing module, \textit{grompp},  which has to be run once before a simulation and thus does not introduce extra computational costs in the simulation itself.

    In the original GROMACS code, the module \textit{grompp} reads the  GROMACS input files and processes them for further use in the molecular dynamics module, \textit{mdrun}. It also checks input data and, if necessary, generates warnings that allow the users to reconsider their chosen setup. For example, the input time step $\Delta t$ is inspected for its ability to provide a stable numerical integration in molecular dynamics. This check is implemented in the \textit{check\_bonds\_timestep($\cdot$)} routine and consists of two main steps. First, for each pair of bonded particles the corresponding period $T$ is calculated with the help of \eqref{eq:harmonicperiod}.
    Then, for given $\Delta t$ and $T$, the Verlet stability condition  $5 \Delta t < T$ \cite{Mazur} is checked. If the condition does not hold, an error message is issued and the simulation is not allowed. (It is easy to see that this restriction is in agreement with condition \eqref{eq:ResonanceStability}, since $1/(\sqrt{2} \pi) \approx 1/5$.) Otherwise, if $10 \Delta t \geq T$ \cite{Mazur} the code issues a message warning that instabilities may arise and recommending to decrease $\Delta t$ or to use a constrained algorithm. Once a warning or error message appears, the search for further problematic oscillations stops.

    For our purposes, we modified this part of the code in such a way that the search continues till the period of the fastest oscillation $\tilde{T}$  is found. Its value is used to define $\bar h$ in \eqref{eq:hbar}.
    Then the optimal parameter value  ${b}$ is  calculated by means of \eqref{eq:Normbar}. A particle swarm optimization algorithm driven by a golden section search \cite{Optimization} is used to perform the required minimization. The parameter ${b}$ is stored in the \textit{input record} structure of GROMACS, so that it can be accessed from every routine in the package after running the \textit{grompp} preprocessing module.

    In standard GROMACS, molecular dynamics simulations are performed with the \textit{mdrun} module using the input file \textit{.tpr} generated by \textit{grompp}.
    The velocity Verlet integrator is implemented in the \textit{update($\cdot$)} function, which is called from \textit{do\_md($\cdot$)} sequentially to update velocities, positions and velocities again. The procedure is repeated as many times as desired.

    To efficiently implement multi-stage integrators in the GROMACS package, it is useful to present a multi-stage scheme in kick/drift factorization form \cite{GROMACSPaper13}.
    For example, two-stage integrators are best rewritten in the form \eqref{eq:rewrite}, which is more suitable for its implementation inside the \textit{mdrun} module in GROMACS.
    The scheme can be implemented with six evaluations of the \textit{update($\cdot$)} function, alternating velocity and position updates in which modified parameters such as $b$, 1/2 and $1/2 - b$ are used.
    Note that this formulation also allows the extension to constrained dynamics, as  explained in \autoref{Ap:Constraints}.
    With our implementation, multi-stage integrators have computational costs equal to those of the standard Verlet method, provided that the latter is run with the choice of time step that equalizes the number of force evaluations.

    The flow chart in \autoref{fig:flowChartAdaptive} summarizes AIA.

\section{Numerical experiments}\label{sec:Testing}
  \subsection{Testing procedure}
    In order to evaluate the efficiency of the proposed AIA scheme, we compared it in accuracy and performance with  the velocity Verlet integrator and with the two-stage integrator (BCSS) of Blanes et al. \cite{BlanesCasasSanzSerna}. In addition, some selected tests also involved the two-stage HOH scheme by Predescu et al. \cite{Predescu}.

    All tests probing various integrating schemes have been repeated with three different simulation techniques, MD, HMC and GHMC. We omit here the data obtained with GHMC for two reasons. First, as expected, HMC and GHMC showed very similar behavioural trends. On the other hand,  the GHMC method possesses an extra parameter that needs to be tuned properly to guarantee  optimal performance. Such a tuning is likely to be time consuming and was not attempted. We therefore decided to avoid reporting  data that may not correspond to the best possible performance of GHMC.

    To ensure a comparison as clear as possible, the following points have been taken into account.

    As we have explained repeatedly (see \autoref{Formulation} for details), whenever a two-stage splitting scheme (AIA or not) and Verlet are used on the same problem,
    the comparisons here are fair (in computational cost terms) because Verlet is run with half the step size and double number of steps.

    In hybrid Monte Carlo (HMC and GHMC) simulations, the number of Metropolis tests was also kept constant regardless of the acceptance rate achieved. For two-stage integrators, the number of MD time steps between two successive Monte Carlo tests was chosen half of the corresponding number for Verlet.

    A broad range of  step sizes has been tested for two benchmark systems with the aim of observing the dependence of the optimal parameter ${b}$ in AIA on  the value of $\Delta t$. Different lengths of MD trajectories in HMC simulations were also explored. Each individual test has been repeated 10 times for unconstrained dynamics and 15 times for constrained dynamics and every single point in the reported data in this paper was obtained by averaging over the multiple runs to reduce statistical errors.

  \subsection{Benchmarks and Simulation setup}
    Two test systems were chosen for the numerical experiments: one describes the non-constrained coarse-grained VSTx1 toxin in a POPC bilayer \cite{Jung2005} and the other the constrained atomistic 35-residue villin headpiece protein subdomain \cite{Villin,NMRVillin}. We will refer to these systems as toxin and villin respectively.

    In the coarse-grained toxin system, four heavy particles on average were represented as one sphere \cite{Wallace2007,Shih2006}, which produced a total number of 7810 particles.
    For both Coulomb and Van der Waals interactions the shift algorithm was used.
    Both potential energies were shifted to 0~kJmol$^{-1}$ at a radius of 1.2~nm.
    Periodic boundary conditions were considered in all directions.
    No constraint algorithm was applied to this system.
    The total length of all simulation runs was 20~ns, which was sufficient, with stable time steps, for a complete equilibration of the system.

    The villin protein was composed of 389 atoms and the system was solvated with 3000 water molecules.
    Coulomb interactions were solved with the PME algorithm of order 6 and Van der Waals interactions were considered as in the toxin system, with the only difference of a radius of 0.8~nm.
    Periodic boundary conditions were again defined in all directions.
    The bonds involving hydrogens were constrained.
    Instead of constraining all  atoms, as it is commonly suggested in the literature (see \cite{vanderSpoelLindahl2003} for instance), we have only constrained the hydrogens, because it is the only case that allows the integration algorithm to perform in parallel with domain decomposition \cite{GROMACSPaper08}.
    Constraining only the hydrogen atoms does not affect the accuracy of the simulation, but  allows  bigger time steps for the integration.
    Since the villin system is an atomistic model, simulations are expected to be slower than for the coarse-grained toxin. However, an exhaustive study of the complete folding process of the villin protein is out of the scope of this work. Thus, with the available computational resources, simulations were run only to observe the effect of the AIA on accuracy and performance of a constrained atomistic system.
    It has to be remarked also that there are examples in the literature of similar tests for which a weak coupling thermostat and a barostat were used to have more realistic results \cite{vanderSpoelLindahl2003}. Barostats are not considered in this study, since the aim is to compare the performance of the AIA scheme with that of  the standard velocity Verlet.
    The total length of all experiments performed for this system was 5~ns.

    The temperature in MD simulations was controlled by the standard v-rescale algorithm for both benchmarks. The reference temperatures were 310~K for toxin and 300~K for villin.
    The same temperatures were used in HMC and GHMC. No thermostat is required in HMC simulations.

\section{Results}\label{Results}
  \begin{figure}[ht]
      \centering
      \subfloat{\includegraphics[width = 3in]{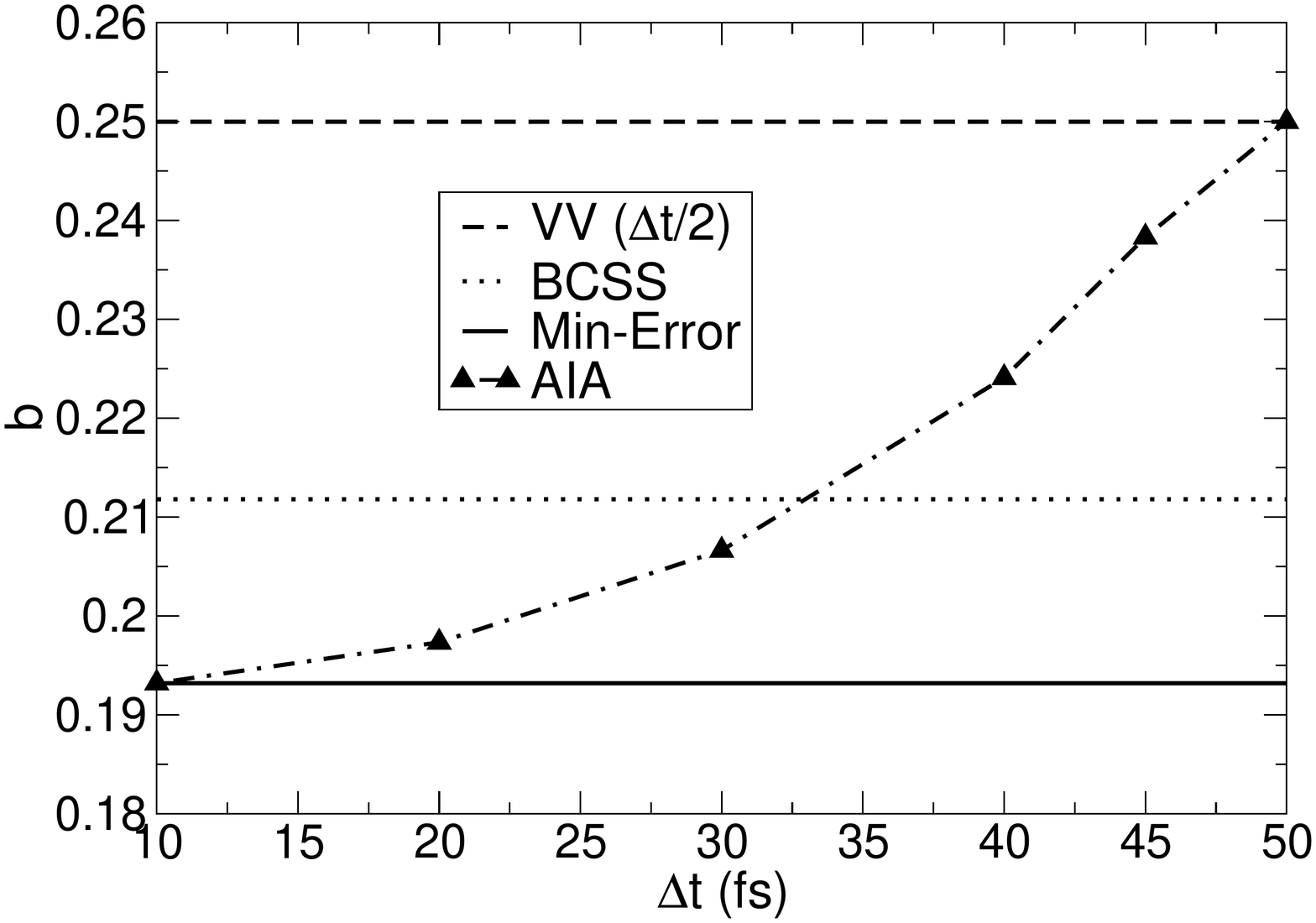}}
      \subfloat{\includegraphics[width = 3in]{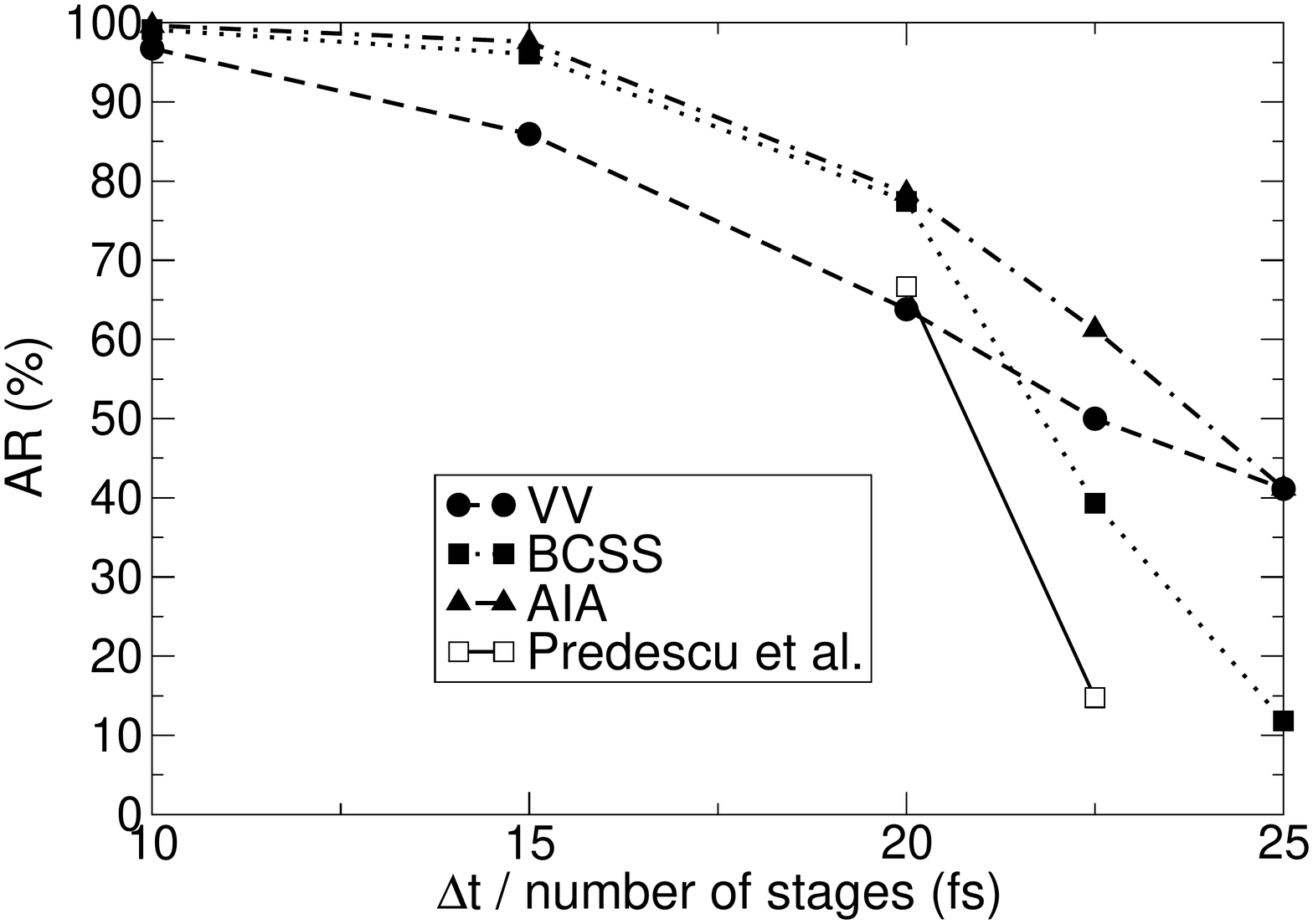}}
      \caption{Dependence of the parameter $b$ on the choice of  $\Delta t$ (left) and its effect on resulting acceptance rates in HMC simulations of the toxin (right).\lq\lq Number of stages" as appears in x-axis label refers to 1 for velocity Verlet and 2 for all two-stage integrators.}\label{plot:bToxin}
    \end{figure}

  We stress that throughout this section the different setups used for the simulation will be expressed in terms of parameters appropriate for the velocity Verlet (one-stage) integrator. This implies that for two-stage schemes the time-steps are doubled and the trajectory lengths are halved which guarantees the fair comparison between these integrators. For improving the readability all the plots have been created following the same criteria.

  \subsection{Unconstrained system}\label{sec:Unconstrained}
    We first present the results for the unconstrained test system.

    The tests were run using the following set of time steps for the Verlet integrator \{10~fs, 15~fs, 20~fs, 22.5~fs, 25~fs\} (recall that for two-stage integrators these values are doubled).
    Two different number of steps in the  MD trajectories, $L$, have been tested in the HMC experiments for each $\Delta t$. In the case of velocity Verlet, the values of $L$ were 2000 and 4000 for all $\Delta t$ except when $\Delta t =25$~fs, where $L$ was chosen to be 1000 and 2000. The corresponding values of $L$ for the two-stage schemes are, as pointed out repeatedly above, halved. The acceptance rates that appear in \autoref{plot:bToxin} were obtained by averaging over all experiments with the same $\Delta t$, regardless of the choice of $L$.

    \begin{figure}[ht]
      \centering
      \subfloat{\includegraphics[width = 3in]{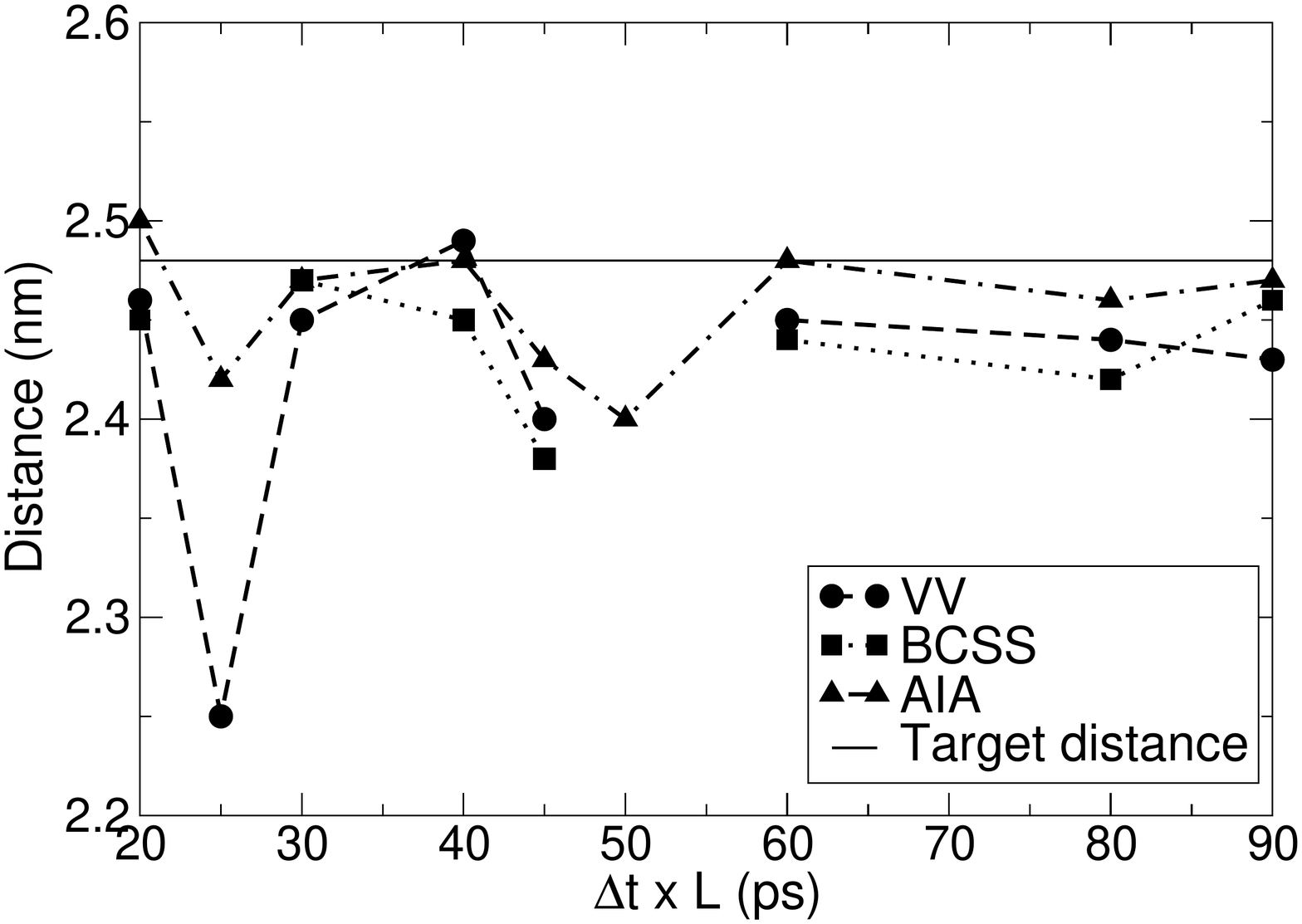}}
      \subfloat{\includegraphics[width = 3in]{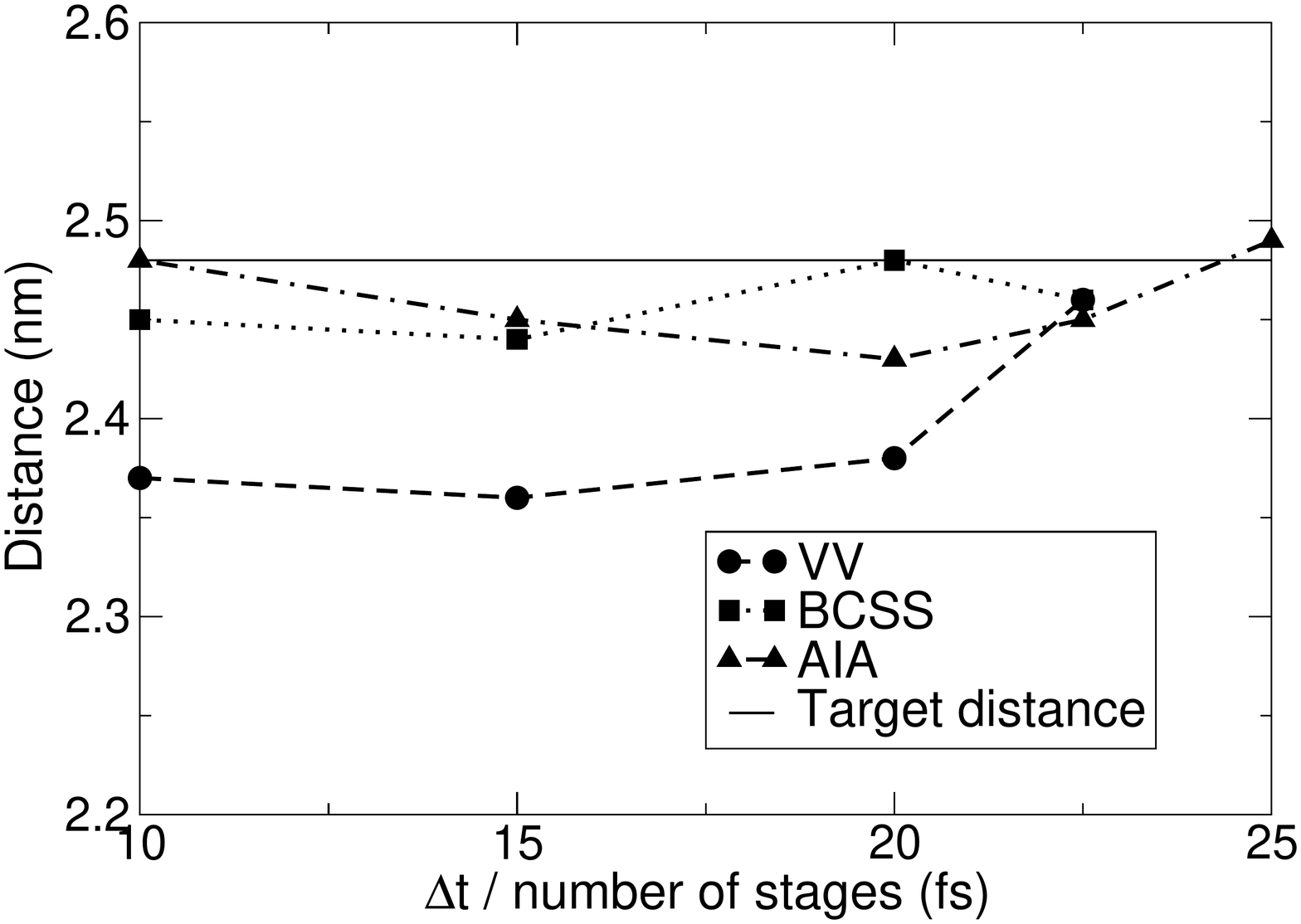}}
      \caption{Distance between the c.o.m. of the toxin and the c.o.m. of the bilayer (expected to be $\sim$2.48 nm) predicted by HMC simulations with different lengths of trajectories $L$, time steps $\Delta t$ and integrating schemes (left) and by MD simulations using various time steps $\Delta t$ and integrators (right).}\label{plot:DistToxin}
    \end{figure}

    As stated earlier,  AIA  finds, for a given physical system and a chosen time step, the unique value of the parameter $b$ in \eqref{eq:IntegratorDef2} that provides the best energy conservation achievable with the members of the family \eqref{eq:IntegratorDef2}.
    \autoref{plot:bToxin} presents the  parameter $b$ determined by the AIA, as a function of $\Delta t$, for simulations of the toxin and compares them with the ones previously identified for different two-stage integrating schemes. As it was intended, for small $\Delta t$, AIA chooses McLachlan's minimum error constant method, and, as $\Delta t$ increases, $b$ approaches $0.25$, a value which, as discussed in section~\ref{Formulation}, essentially yields the Verlet integrator.
    The two-stage integrator BCSS  \cite{BlanesCasasSanzSerna} is the optimal choice for time steps roughly twice smaller than the stability limit of the velocity Verlet integrator.

    We then investigated the effect of the AIA  on the performance of HMC simulations by monitoring acceptance rates as functions of $\Delta t$  with different two-stage integrators.
    Conservation of energy has a direct impact on acceptance or rejection in the Metropolis test of the hybrid Monte Carlo methods: the better the energy is preserved, the more proposed trajectories are accepted \cite{BeskosPillaiRobertsSanzSernaStuart}. Thus, by design,  AIA  has to provide, at least for Gaussian distributions, the highest acceptance rates for any choice of $\Delta t$. This is demonstrated in \autoref{plot:bToxin}. The two-stage schemes of  \cite{BlanesCasasSanzSerna} and \cite{Predescu} ensure higher acceptance rates than velocity Verlet for  time steps significantly smaller than the Verlet stability limit. However the performance of those two-stage schemes drops  significantly for larger time steps.  AIA  yields acceptance rates that are as good as those of BCSS  when $\Delta t$ is small and as good as those of Verlet near the Verlet stability limit. In particular AIA does not yield worse results than Verlet for any values of $\Delta t$.

    The trend observed for the HMC method as shown in \autoref{plot:bToxin} was also apparent in GHMC tests.

    \begin{figure}[ht]
      \centering
      \subfloat{\includegraphics[width = 3in]{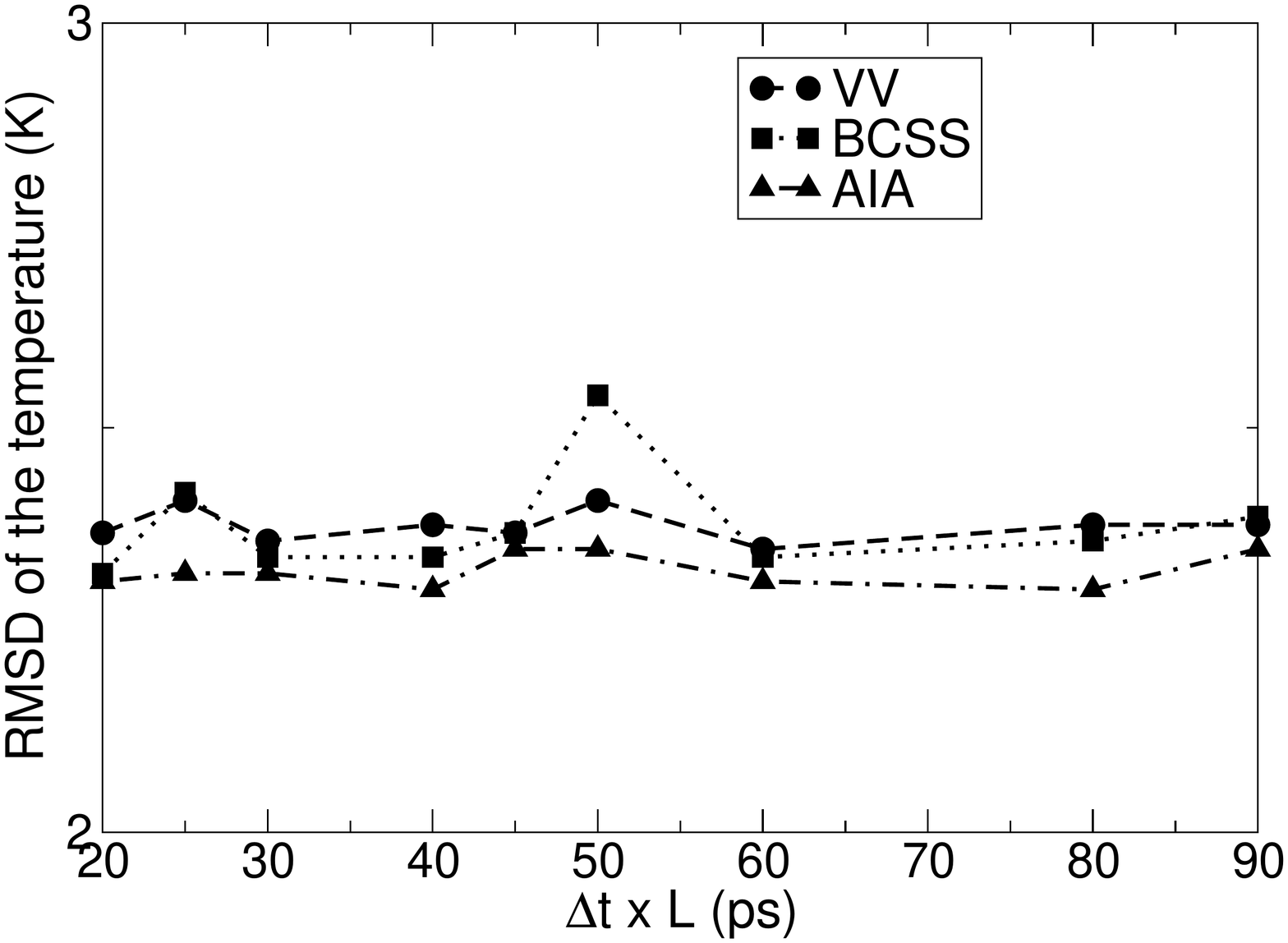}}
      \subfloat{\includegraphics[width = 3in]{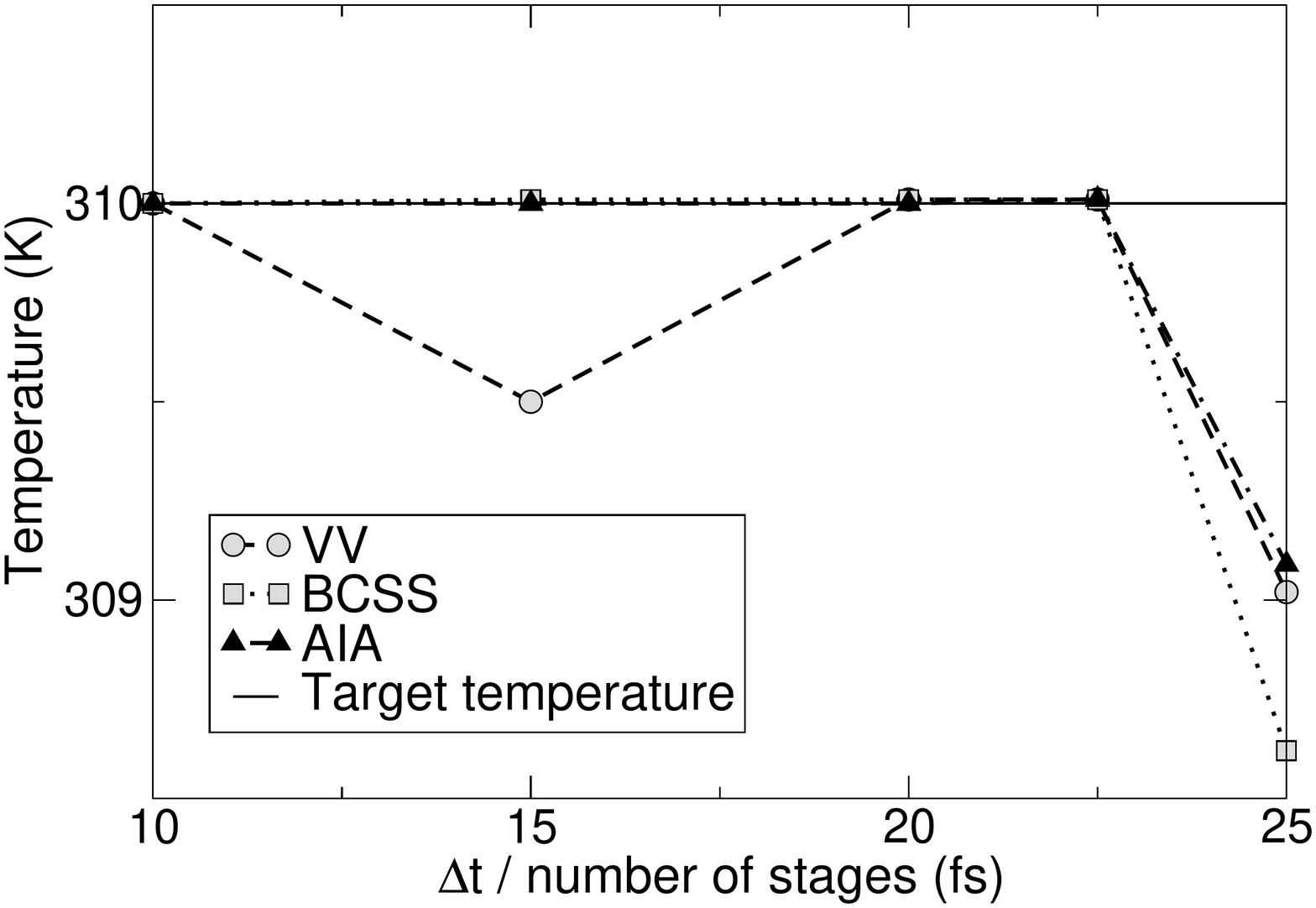}}
      \caption{Temperature RMSD with respect to the target temperature observed in HMC simulations of the toxin with different lengths of trajectories $L$, time steps $\Delta t$ and integrating schemes (left) and average temperature in MD simulations of the toxin using various time steps $\Delta t$ and integrators (right). The target temperature was set to 310~K. The v-rescale thermostat was applied in MD.}\label{plot:TempToxin}
    \end{figure}

    To compare the impact of different integrating schemes on the accuracy of HMC and MD simulations, we calculated averages for two thermodynamic observables: the temperature $T$ and the distance $d$ traveled by the toxin from the center of the membrane to the preferable location at the surface of the membrane. The expected average values of the distance are around $\sim$2.48 nm \cite{Jung2005,Wee2008}, whereas the target temperature was chosen to be 310~K. The  performed simulations had a fixed total length of 20~ns, which was long enough for equilibrating the system if stable time steps were used, but not sufficient for obtaining accurate averages. So, the tests are meaningful for observing trends rather than obtaining  good production results. For HMC we found more informative to plot  the RMSD between the target temperature and the observed temperatures  rather than the average temperatures themselves.  For MD simulations, where the overall fluctuations are smaller and the trends for averages, even in  short simulations, are clearer, we plot temperatures.

    \autoref{plot:DistToxin} and \autoref{plot:TempToxin} summarize the averages for the two observables, distance and temperature. From now on, we  plot the properties obtained with HMC simulations versus the product $\Delta t \times L$ of the time step and the number of  steps in an MD trajectory. This is due to the important role this product plays in the overall acceptance rate and in the correlation in HMC simulations \cite{nawaf}.

    \begin{figure}[ht]
     \centering
     \subfloat{\includegraphics[width = 3in]{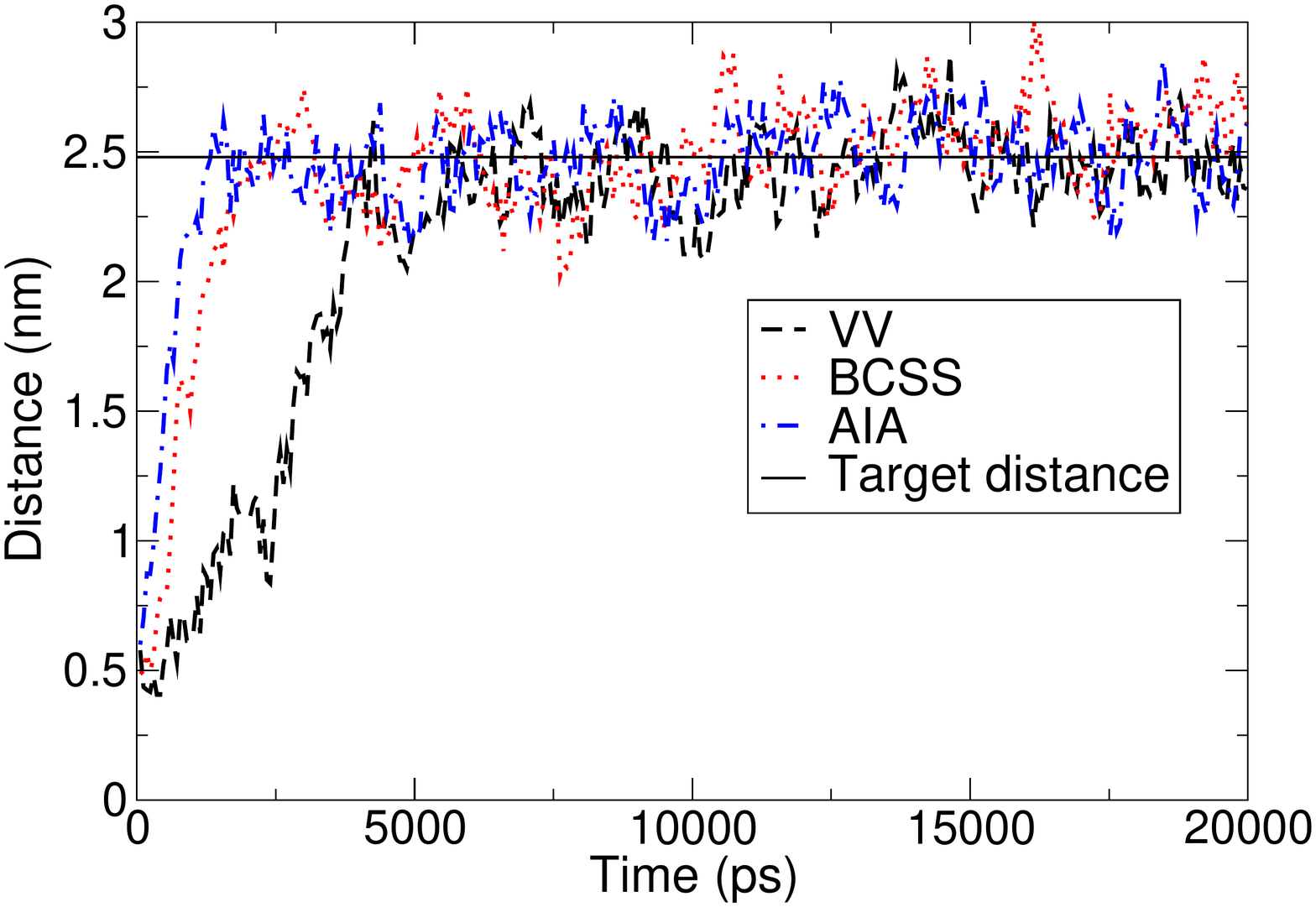}}
      \subfloat{\includegraphics[width = 3in]{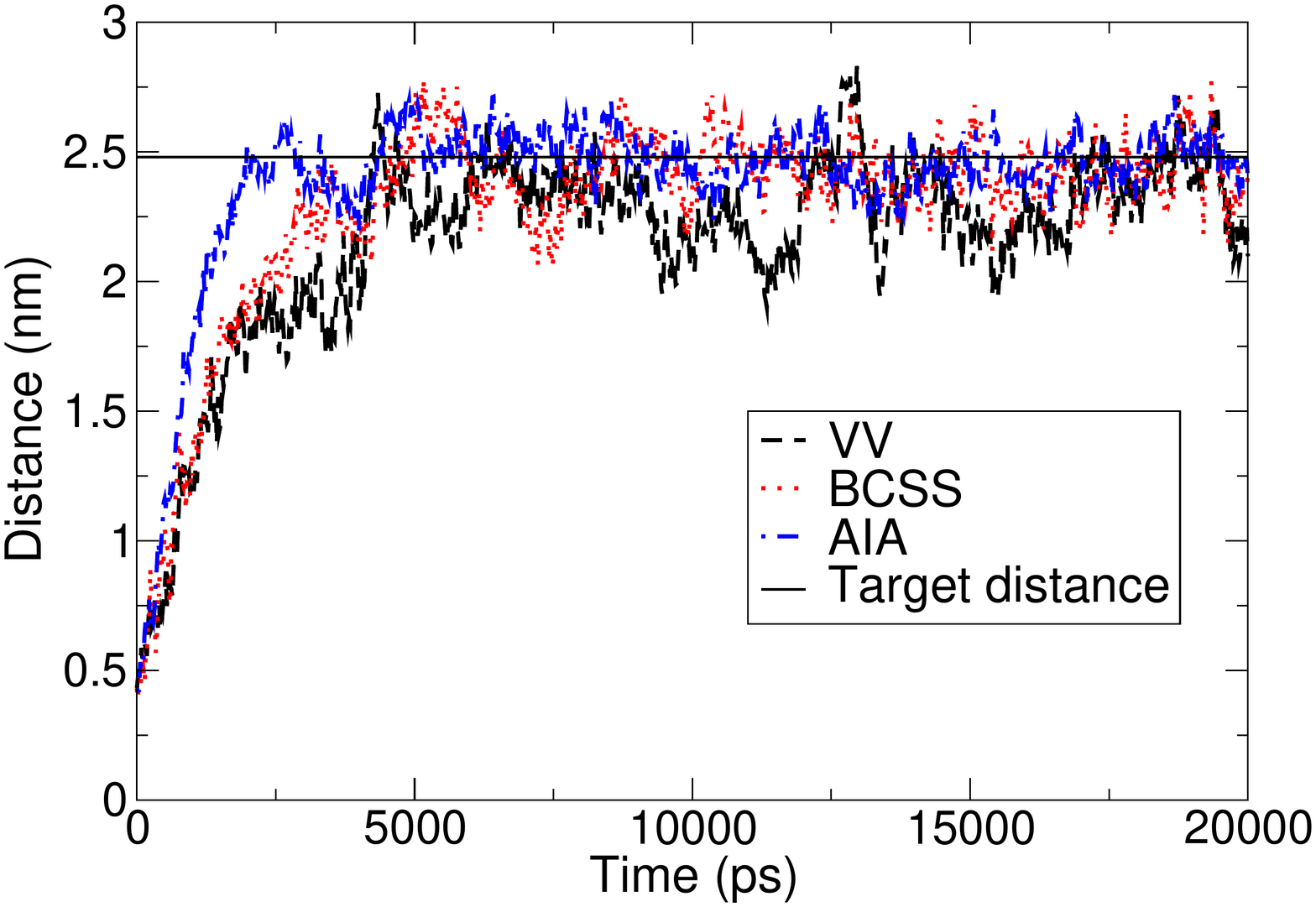}}
      \caption{Distance between the c.o.m. of the toxin and the c.o.m. of the bilayer as a function of time obtained in HMC simulations with time step $\Delta t = 15$~fs, trajectory length $L = 4000$ and different integrators (left) and in MD simulations with time step $\Delta t = 10$~fs and the same integrators (right). The expected value is $\sim$2.48~nm.}\label{plot:DistanceToxin}
    \end{figure}

    As follows from \autoref{plot:DistToxin} and \autoref{plot:TempToxin}, for both  properties, $d$ and $T$, the accuracy of AIA is comparable to, but typically better than, the accuracy provided by BCSS and Velocity Verlet for step sizes distant from the Verlet stability limit. However near the stability limit the accuracy of all integrators decreases - more dramatically for BCSS and less noticeably for AIA. Interestingly, longer MD trajectories ($L = 2000$) in HMC allow AIA to be accurate at such large values of $\Delta t$  (see \autoref{plot:DistToxin} at $\Delta t \times L=50$~ps).
    In contrast, the accuracy in simulations with BCSS and Verlet  are rather sensitive to the choice of $\Delta t$. The former failed to produce meaningful averages for $\Delta t = 25$~fs. Less dramatic differences but  similar trends were observed for molecular dynamics simulations (right panels of \autoref{plot:DistToxin} and \autoref{plot:TempToxin}).

    Finally, we inspected the role of numerical integrators in the sampling efficiency of HMC and MD simulations.

    In \autoref{plot:DistanceToxin} the distance $d$ between the c.o.m. of the toxin and the c.o.m. of the bilayer is shown as a function of time for a single choice of the time step $\Delta t = 15$~fs and the trajectory length $L = 4000$ in HMC, and for $\Delta t = 10$~fs in MD. The superiority of the AIA method is clearly demonstrated in both HMC and MD, since AIA makes the toxin reach the target destination earlier than the rest of the integration schemes do.

\begin{figure}[ht]
      \centering
      \includegraphics[width = 3in]{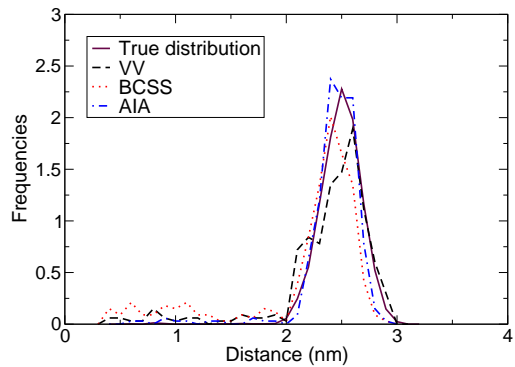}
      \caption{Distribution of the distances between the c.o.m. of the toxin and the c.o.m. of the bilayer observed in HMC simulations of 20 ns length with time step $\Delta t = 15$~fs, trajectory length $L = 4000$ using different integrators. The solid magenta line presents the\lq\lq true\rq\rq distribution produced with a ten times longer simulation (200 ns) that used the same input. The $y$-axis presents frequencies which are calculated as the normalized numbers of hits registered for a distance bin within a simulation. Here normalization is performed with respect to a product of a total number of samples and the size of a distance bin (0.1 in this particular case).}\label{plot:Distribution}
    \end{figure}

    \autoref{plot:Distribution} presents the distributions of the distances $d$ collected from  simulations with different integrators, i.e. AIA, VV and BCSS, and compares them with the \lq\lq true\rq\rq distribution obtained from the HMC simulation at $\Delta t = 15$~fs and $L = 4000$ of 200 ns length, i.e. ten times longer than the other ones. It can be seen that AIA samples more closely to this distribution. As for all tests in this section, the plotted data are resulted from averaging over several repetitive runs (see \autoref{sec:Testing} for more details).

    Finally, the integrated autocorrelation function IACF of the drift of the toxin to the preferred interfacial location was measured during the equilibration stage of the simulations {for the range of step sizes and trajectory lengths. The autocorrelation function (ACF) is a commonly used tool for evaluating sampling efficiency in molecular dynamics simulations \cite{AllenTildesley}, statistics and other fields. For a certain property $f$ depending on time it is defined as
    $$\mbox{ACF}(f(t)) = \langle f(\xi) f(\xi + t) \rangle_{\xi}.$$
    The integral of the correlation function over time is called the integrated autocorrelation function (IACF) or integrated autocorrelation time
    $$\mbox{IACF}(f(t)) = \int_0^{\infty} \mbox{ACF}(f(t)) \mathrm{d}t.$$
    Intuitively, IACF can be understood as measuring the time needed, on average, for generating a non correlated sample. It can be seen as the inverse of the effective sample size (ESS) \cite{ESS},  a measure often used in statistical applications of Monte Carlo methods. In practice, all the correlation functions are calculated for discrete values. Low values of measured IACFs mean low correlations between the generated samples and thus better sampling.

    \begin{figure}[ht]
      \centering
      \subfloat{\includegraphics[width = 3in]{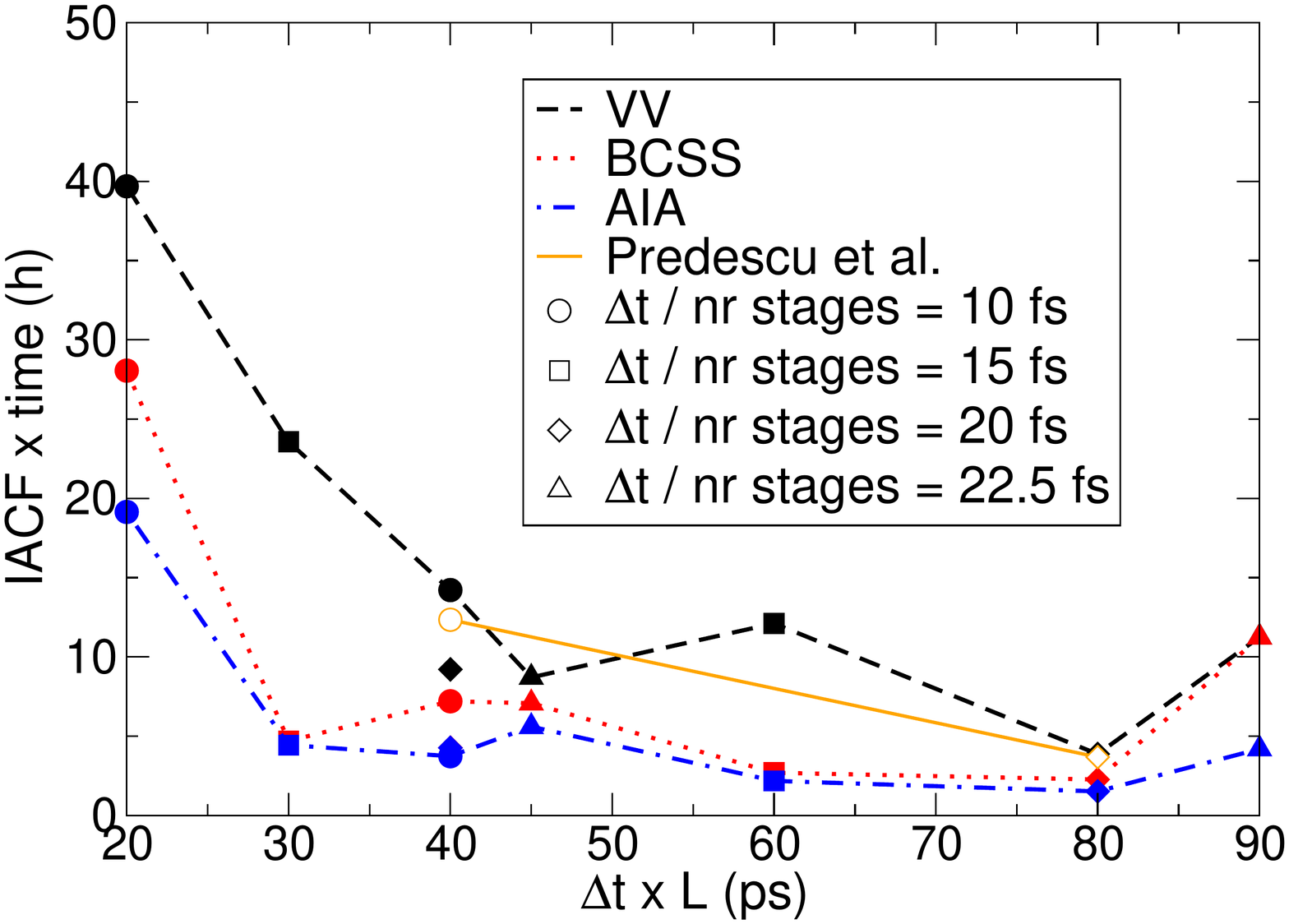}}
      \subfloat{\includegraphics[width = 3in]{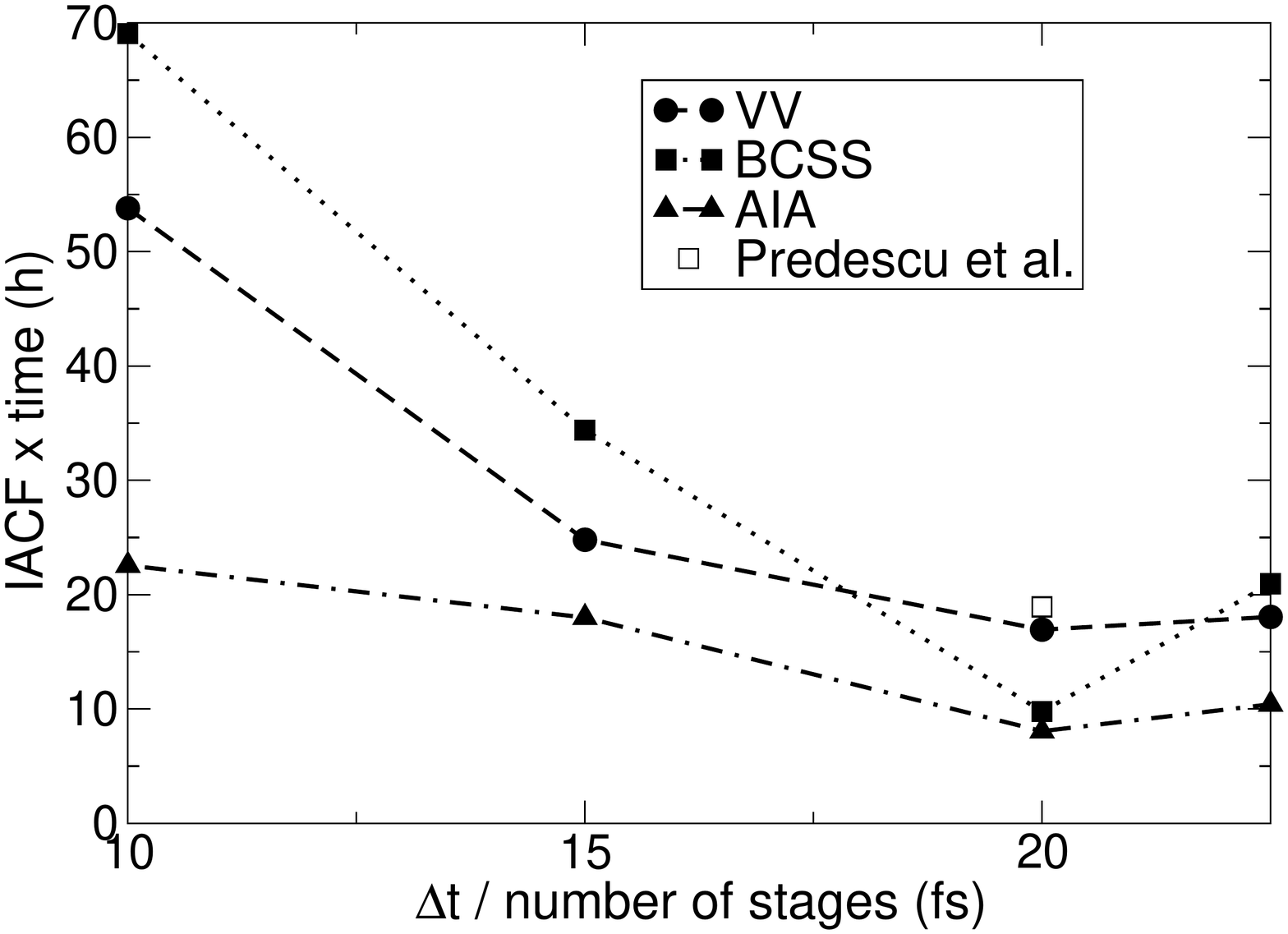}}
      \caption{IACF of the drift of the toxin to the preferred interfacial location evaluated as a function of  $L$ and $\Delta t$ in HMC tests (left) and as a function of $\Delta t$ in MD runs (right). Four integrating schemes were tested in HMC and MD simulations: velocity Verlet (a dashed line), the two-stage integrator BCSS  (a dotted line), the HOH-integrator of Predescu et al. (open symbols) and the AIA integrators (a dot-dashed line).}\label{plot:IACFequil}
    \end{figure}

    \autoref{plot:IACFequil} presents the IACF measured in HMC and MD with different integrating schemes for the same range of time steps and trajectory lengths  described above. Note, in the vertical axis, that computing time is used to normalize the results. The IACF values for 25~fs/50~fs are not plotted since the lack of stability at those  step lengths in all integrating schemes  produces poor, non-informative results. For completeness we present the data in the Appendix.

    In \autoref{plot:IACFequil} we use different symbols for different values of $\Delta t$  to provide a better feeling for the relation between $\Delta t$ and the efficiency achieved. Two different symbols corresponding to the same  $\Delta t \times L$ mean that two different combinations of $\Delta t $ and $L$ are possible to get the same number on the $x$ axis.

    As seen from  \autoref{plot:IACFequil}, for all combinations of $\Delta t$ and  $L$, both HMC and MD simulations using the AIA integrators decorrelated faster  than the corresponding simulations that used the velocity Verlet integrator, BCSS or the method of Predescu et al. In fact, for some specific choices of $\Delta t$ the AIA integrators led to an efficiency several times higher than that of the velocity Verlet or any of the tested two-stage integrators. This applies to both simulation methods, HMC and MD. The fact that the better energy conservation of  AIA led to  better sampling efficiency in hybrid Monte Carlo simulations was  not surprising. For molecular dynamics,  better conservation energy guarantees  better accuracy but not necessarily better sampling. However \autoref{plot:IACFequil} clearly demonstrates the positive impact of  energy conservation on the sampling performance of MD. Still, comparison of the two plots in \autoref{plot:IACFequil} reveals the clear superiority  in sampling efficiency of HMC over MD for the tested system.

    \begin{figure}[ht]
      \centering
      \subfloat{\includegraphics[width = 3in]{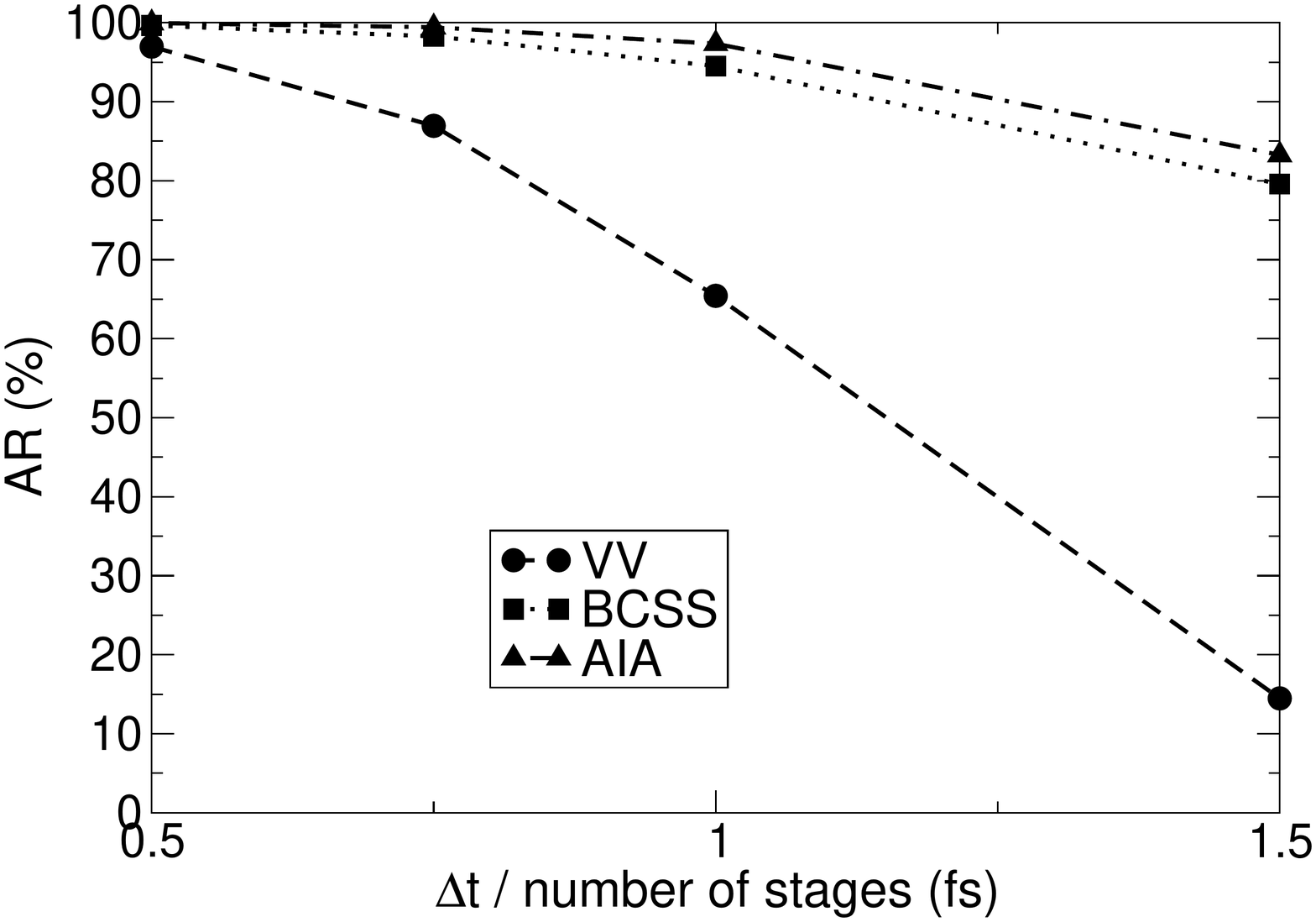}}
      \subfloat{\includegraphics[width = 3in]{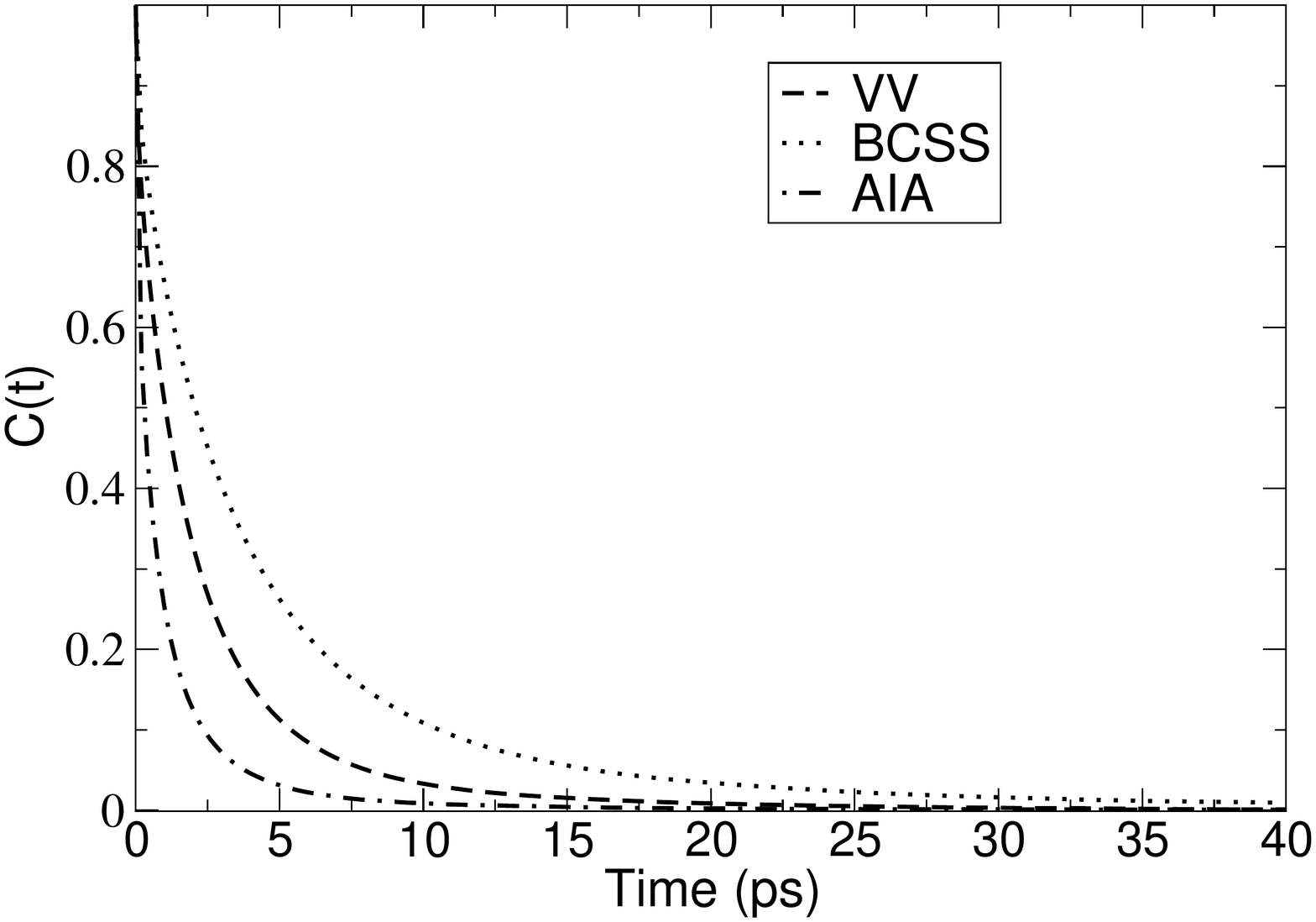}}
      \caption{Effect of the paramter $b$ on the resulting acceptance rates in HMC simulations of water (left) and autocorrelation functions of the hydrogen bonding in MD simulations (right) for $\Delta t = 2$~fs. The two-stage integrator looses performance at the chosen time step whereas the AIA not only outperforms this integrator but also shows faster convergence than the standard velocity Verlet provides. The IACF's are: VV = 12.31, BCSS = 22.92, AIA = 5.66.}\label{plot:Water}
    \end{figure}

    A few more  useful observations may be extracted from \autoref{plot:IACFequil}. Analyzing the IACF calculated for HMC simulations with different combinations of $\Delta t$ and $L$, one can conclude that, for fixed $\Delta t$, a larger $L$ gives better performance for all integrators. Moreover, to achieve better performance, the choice of the product of $\Delta t$ and $L$ is more important than $\Delta t$ itself. For instance, $\Delta t = 30$~fs and $L = 2000$ is a better choice than $\Delta t = 40$~fs and $L = 1000$.

    At this stage, we can conclude that the Adaptive Integration Approach outperforms the other tested schemes in accuracy, stability and sampling efficiency for all tested step sizes. As one can expect, long  step sizes, close to the maximum allowed by stability, lead to  accuracy and performance degradation in all schemes. For the adaptive scheme this effect is much smoother.

    These conclusions are also supported by the results obtained in HMC and MD simulations of 216 molecules of water at 300 K. The model used is the flexible version of SPC \cite{SPC}.
    Taking into account the important role water plays in bimolecular simulations, we include here two plots in \autoref{plot:Water} showing the advantage of  AIA  over other integrating schemes in sampling with HMC (left) and MD (right) simulations.

  \subsection{Constrained system}
    \begin{figure}[ht]
      \centering
      \subfloat{\includegraphics[width = 3in]{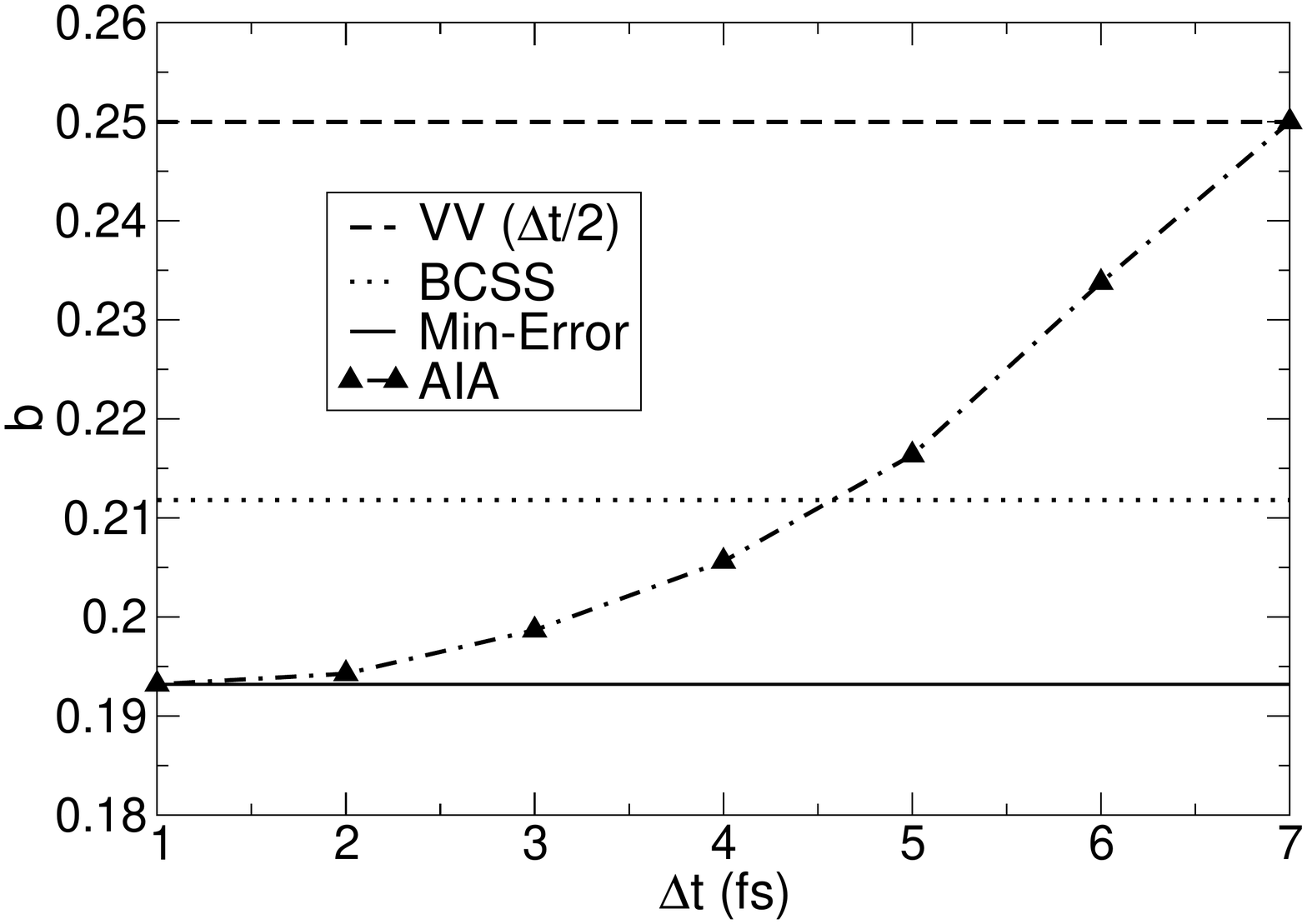}}
      \subfloat{\includegraphics[width = 3in]{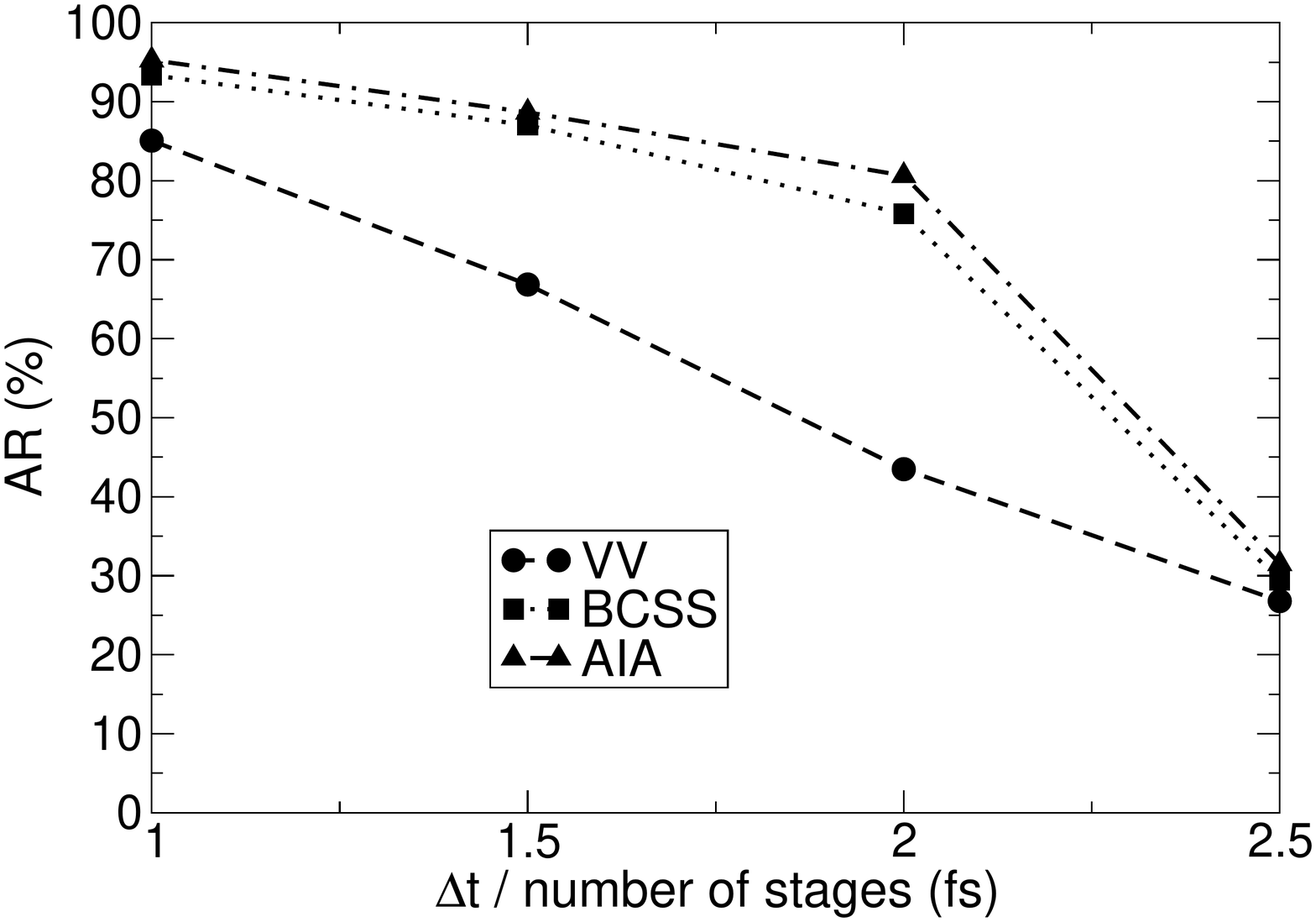}}
      \caption{Dependence of the parameter $b$ on the choice  $\Delta t$ (left) and its effect on the resulting acceptance rates in HMC simulations of villin (right).}\label{plot:bVillin}
    \end{figure}

    For testing efficiency of the AIA integrators in simulations of constrained systems we followed the same strategy as in \autoref{sec:Unconstrained}. The time steps chosen for the tests in this case, however, were in the range typical for step sizes used in atomistic simulations and thus differed from those considered in coarse-grained experiments in \autoref{sec:Unconstrained}. More specifically, we tested the following time steps, $\Delta t / nr$ ($nr$=1 for Verlet and 2 otherwise): 1~fs, 1.5~fs, 2~fs, 2.5~fs. The number $L$ of  steps in MD trajectories in HMC were exactly the same as in \autoref{sec:Unconstrained}, i.e.\ 2000 and 4000 in the tests with Verlet, and 1000 and 2000 for the two-stage methods. The measured acceptance rates were averaged over different lengths $L$ for each $\Delta t$.

    To our satisfaction, the  positive impact of the AIA strategy on the quality of simulations  demonstrated in unconstrained systems has  also been observed in the case of constrained dynamics.

    \begin{figure}[ht]
      \centering
      \subfloat{\includegraphics[width = 3in]{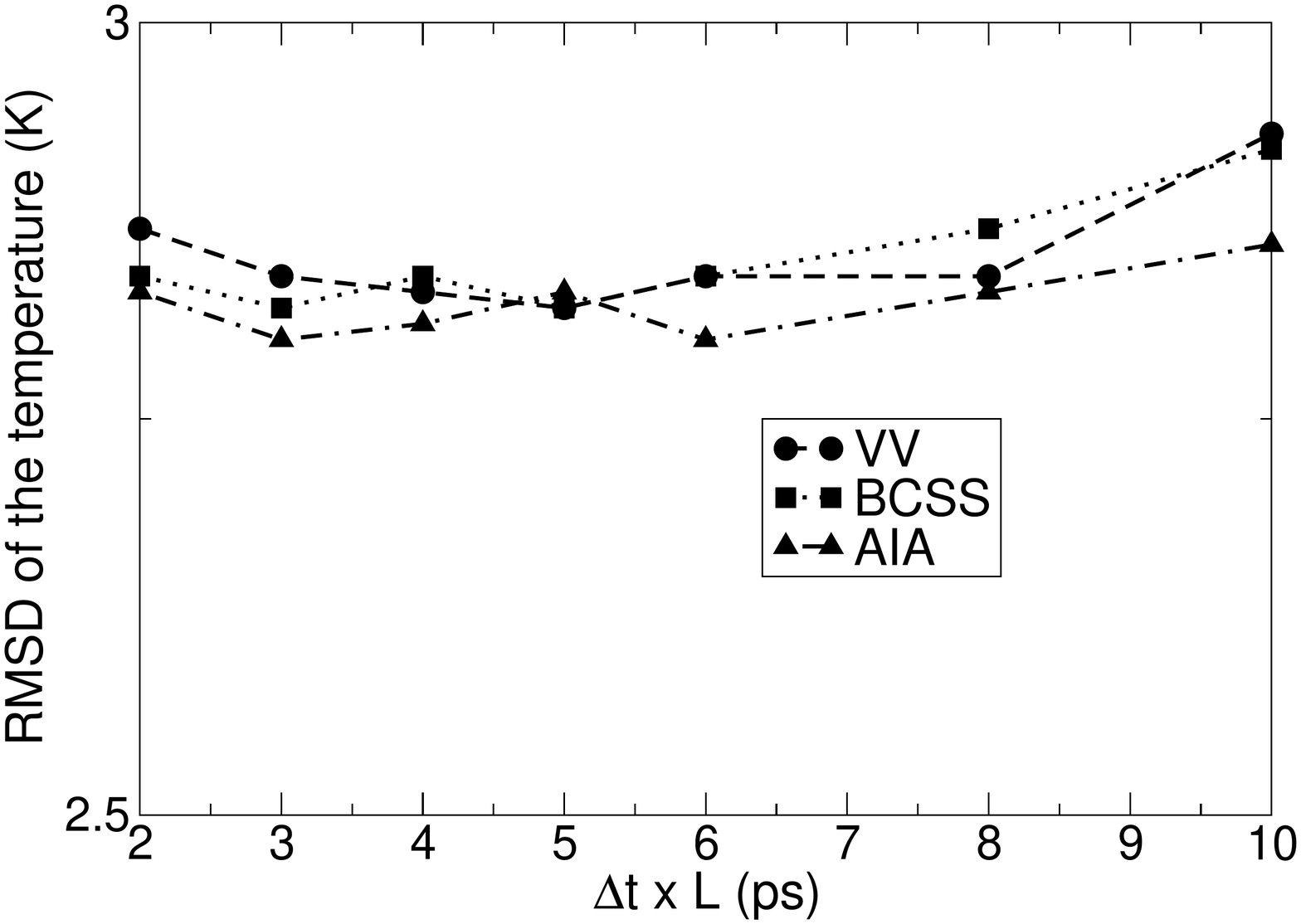}}
      \subfloat{\includegraphics[width = 3in]{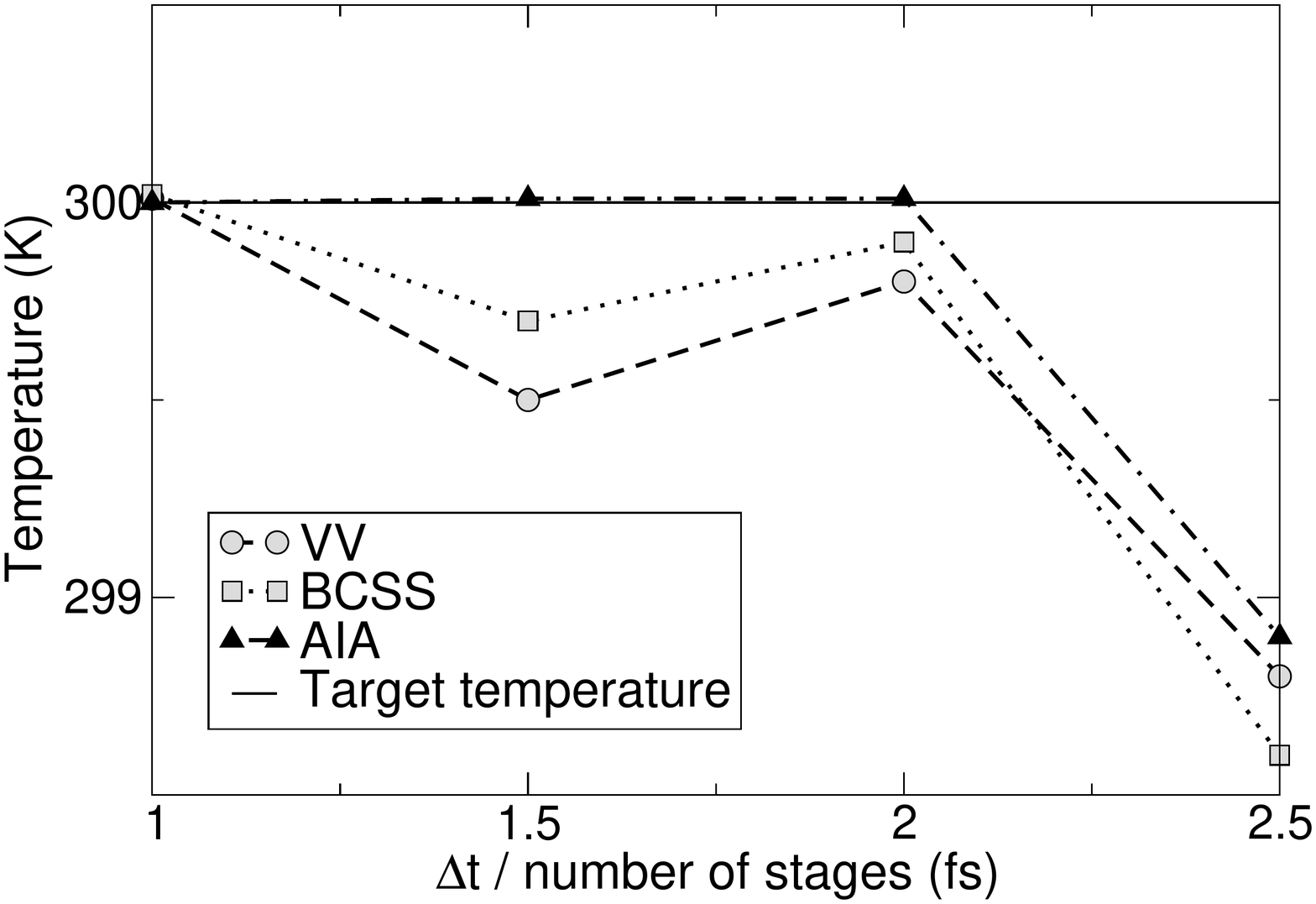}}
      \caption{Temperature RMSD with respect to the target temperature in HMC simulations of villin with different lengths of trajectories $L$, time steps $\Delta t$ and integrating schemes (left) and average temperature in MD simulations of villin using various time steps $\Delta t$ and integrators (right). The target temperature was set to 300~K. The v-rescale thermostat was applied in MD.}\label{plot:TempVillin}
    \end{figure}

    \autoref{plot:bVillin} shows  trends that match  those summarized in \autoref{plot:bToxin}. The only significant difference is  for BCSS; where  the  loss in performance at larger $\Delta t$ is smaller for villin (\autoref{plot:bVillin}) than for toxin (\autoref{plot:bToxin}).

    As in the case of the unconstrained system, the \lq\lq convergence" of  AIA to the velocity Verlet integrator was also observed (at around 2.25~fs / 4.5~fs), but the resulting acceptance rates were so low in all tests that the corresponding experiments have been excluded from consideration.

    The villin system is a popular benchmark for studying folding processes, due to its comparatively fast folding times. In this paper, we did not aim to investigate  in full the folding of villin. Rather, the fast folding helped us to design  computationally feasible tests for measuring accuracy and efficiency of the different numerical integrators.

    Calculated averages of simulated temperatures in HMC and MD tests were used for evaluating the accuracy provided by the Velocity Verlet integrator and the two-stage integrating schemes of interest. As in \autoref{sec:Unconstrained}, the length of tests with HMC and MD simulations was fixed and sufficient to analyse the effect of $\Delta t$ on the level of accuracy achieved in simulations, but not to guarantee low statistical errors.

    As in \autoref{sec:Unconstrained}, \autoref{plot:TempVillin} shows the dependence of the temperature RMSD with respect to the target temperature on the chosen $\Delta t$, trajectory lengths and integrators for HMC and the average temperatures with the v-rescale thermostat for different time steps and integrators in MD. Evidently, AIA  provided the smallest fluctuations of averages as a function of $\Delta t$ within the inspected range of step sizes, even though the differences in the data obtained with the different integrators were less marked than in the case of the toxin in \autoref{sec:Unconstrained}. Degradation of accuracy was observed for larger $\Delta t$ in all simulations, but  was less visible for  AIA  than for  Verlet or BCSS. The data collected at $\Delta t / nr = 2.5$~fs showed poor accuracy for all tests.

    \begin{figure}[ht]
      \centering
      \subfloat{\includegraphics[width = 3in]{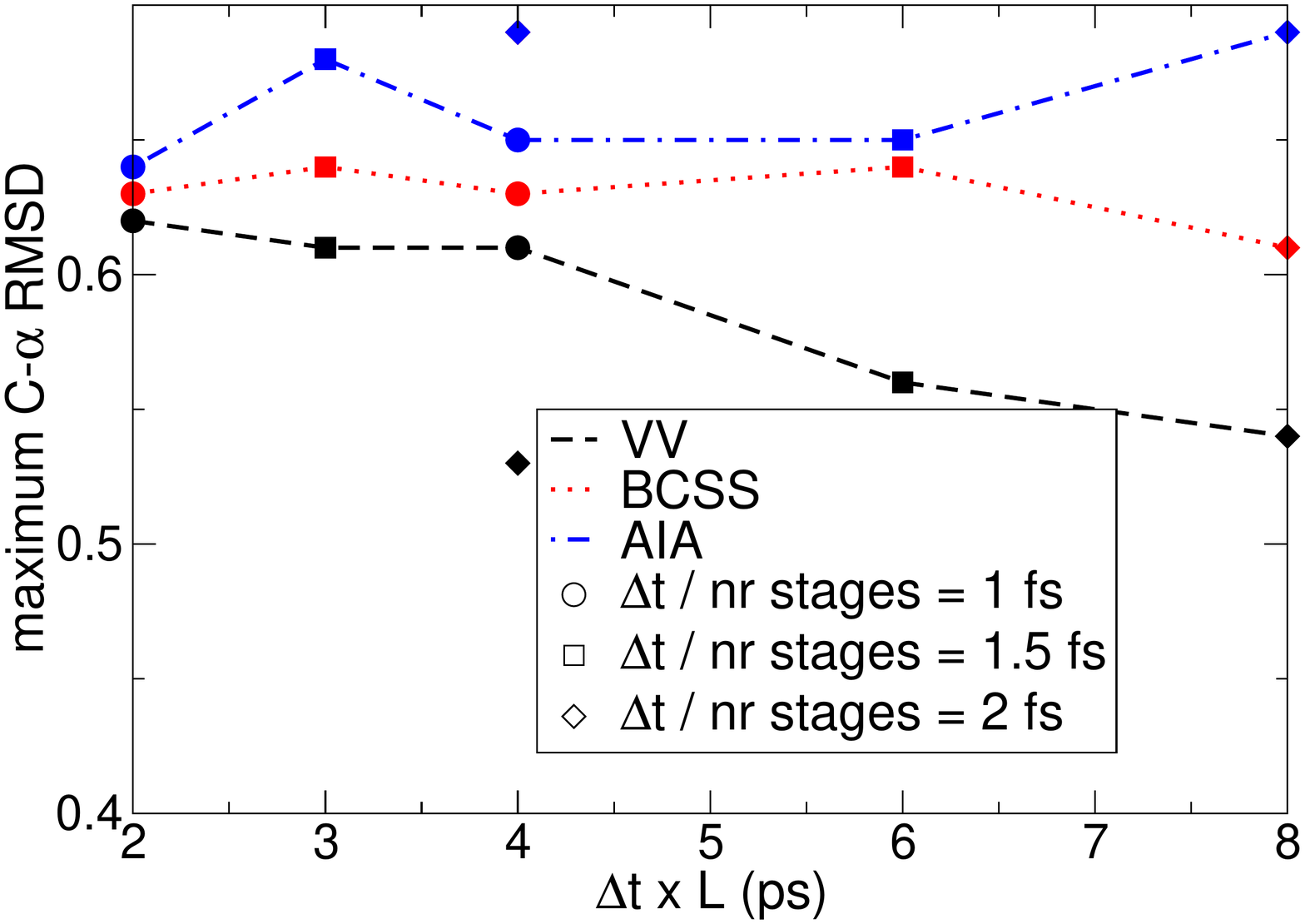}}
      \subfloat{\includegraphics[width = 3in]{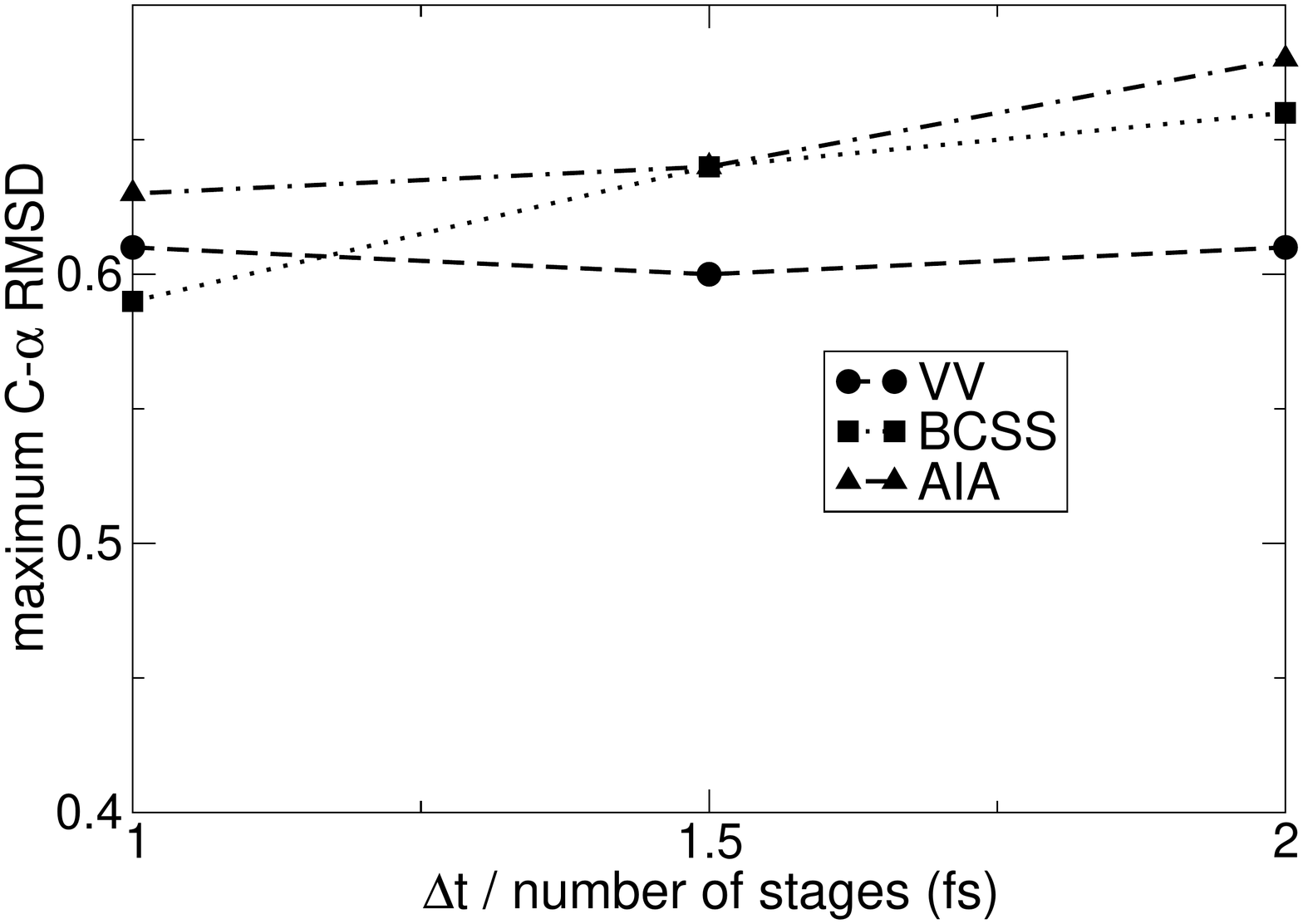}}
      \caption{Maximum $\alpha$-carbon RMSD between any two structures in HMC simulations of villin with different lengths of trajectories $L$, step sizes $\Delta t$ and integrating schemes (left) and in MD simulations of villin using various step sizes $\Delta t$ and integrators (right).}\label{plot:RMST}
    \end{figure}

    We completed our testing of AIA for constrained dynamics with an analysis of its impact on the sampling performance of HMC and MD.
    We chose to measure the quality of sampling through the positional RMSD from the native structure as a function of the simulation steps in both HMC and MD cases.
    The state of a protein folding can be understood by computing the root-mean-square deviation (RMSD) of the $\alpha$-carbon.
    It can be used to make a comparison between the structure of a partially folded protein and the structure of the native state. The RMSD of certain atoms in a molecule with respect to a reference structure is calculated as
    $$\mbox{RMSD} = \sqrt{\frac{1}{N} \sum_{i = 1}^N \delta_i^2},$$
    where $\delta_i$ is the distance between the atoms $i$ in the two structures compared. As it is done in \cite{vanderSpoelLindahl2003}, we have calculated what the authors call RMST, the maximum RMSD of the $\alpha$-carbon between any two structures in a simulation. The idea is to roughly measure the extent of the conformational space sampled in a simulation. As in the unconstrained case, we have also plotted these values for the different combinations of time step and length of trajectories $\Delta t \times L$. In \autoref{plot:RMST} the simulation results obtained with different integrators are compared. It can be observed, in both HMC and MD cases, that AIA leads to broader sampling of the conformational space no matter the choice of time step or trajectory length. The largest difference with respect to velocity Verlet can be observed when the biggest time step $\Delta t = 2$~fs is used.

    We have also computed the radius of gyration, which provides an estimation of the compactness of a desired structure. As in \cite{vanderSpoelLindahl2003}, we have considered the experimental value 0.94~nm \cite{NMRVillin} as a target value. The simulations performed are not long enough to observe any proper convergence to the value, however the tendency of the protein evolution can be seen  through the comparison of the simulated radius of gyration with the target one. In \autoref{plot:gyrate} the average radii of gyration obtained from HMC (left) and MD (right) simulations using different integrators and different values of simulation step sizes and trajectory lengths are presented. While the results associated with the Velocity Verlet and BCSS integrators are still far from the target value, the averages produced with AIA are, regardless a choice of simulation parameters, always closer to 0.94~nm both in HMC and in MD.

    \begin{figure}[ht]
      \centering
      \subfloat{\includegraphics[width = 3in]{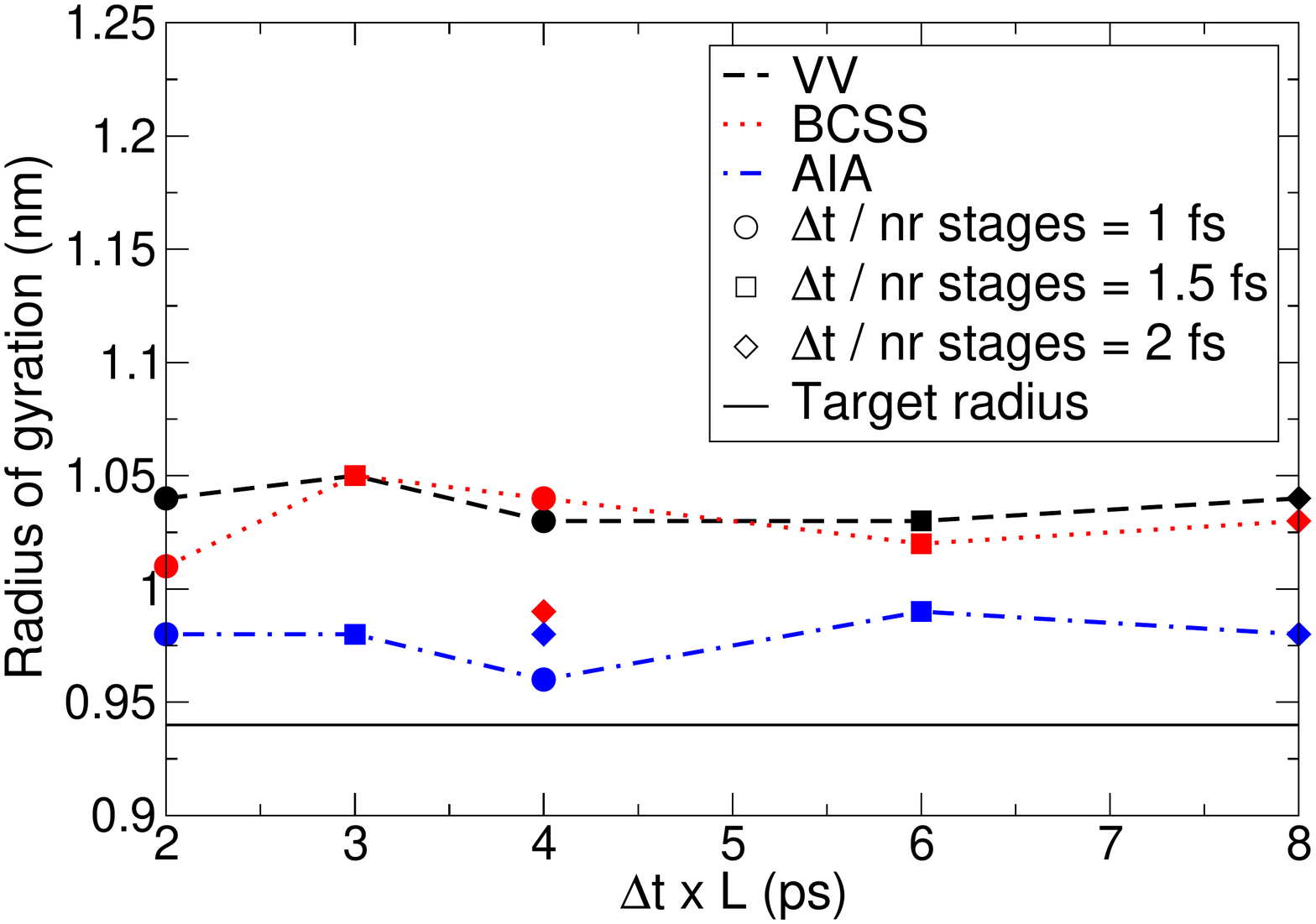}}
      \subfloat{\includegraphics[width = 3in]{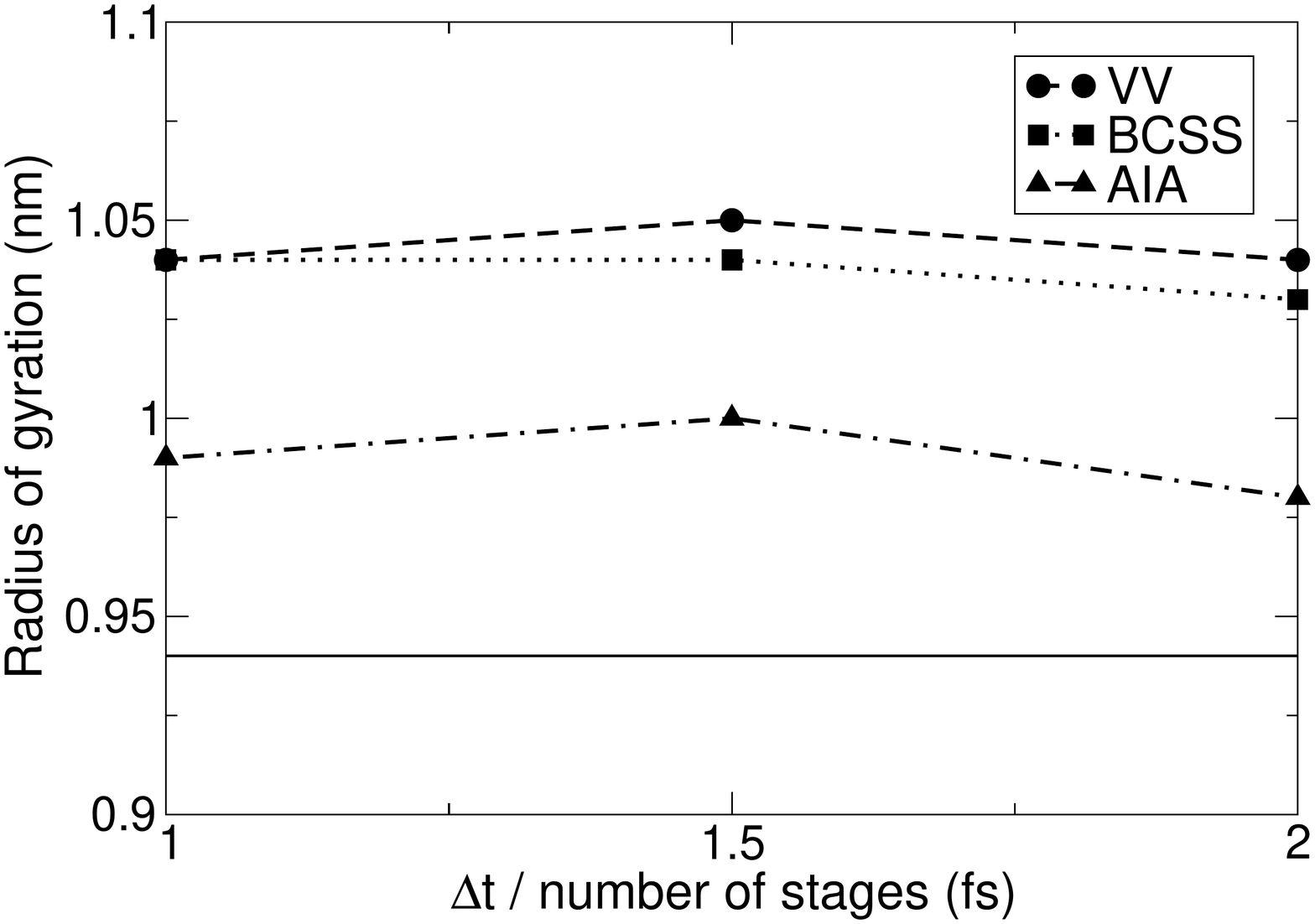}}
      \caption{Average radii of gyration in HMC simulations of villin with different time steps $\Delta t$, lengths of trajectories $L$ and integrating schemes (left) and in MD simulations of villin using various step sizes $\Delta t$ and integrators (right). The target experimental radius of gyration is 0.94~nm.}\label{plot:gyrate}
    \end{figure}

    Similar trends were seen in GHMC simulations. The results are not shown (see \autoref{sec:Testing}).

    Obviously, it is impossible,  with basis on these short tests, to make precise conclusions about  features of the folding process, e.g. about the folding rate.  More detailed studies of the protein folding are advisable. However, what can be concluded without hesitation is that sampling in molecular simulations of atomistic constrained systems with HMC and MD benefits from integrators that guarantee the best possible conservation of energy, as is the case with AIA.

\section{Conclusions}\label{Conclusions}
  We have presented an alternative to the standard velocity Verlet integrator, known to be the state-of-the-art method for numerical integration of the Hamiltonian equations in molecular dynamics. The novel methodology, which we call the Adaptive Integration Approach, or AIA, offers, for any chosen step size, a system-specific integrator which guarantees the best energy conservation for harmonic forces achievable by an integrator from the family of two-stage splitting schemes, including Verlet. While improvements in energy   conservation do not necessarily imply dramatic changes in sampling, they improve acceptance rates in hybrid Monte Carlo methods. The experiments performed in the present study also show that in molecular dynamics AIA leads to improvements of sampling as measured by the metrics considered. The improved sampling may arise as a consequence of either  enhanced accuracy with a given step size or to the possibility of longer step sizes.
  
  The AIA scheme can be implemented, without introducing computational overheads in simulations, in any software package which includes MD and / or HMC. In this study, we implemented the AIA method in multiHMC-GROMACS, a modified version of the popular GROMACS code,  and tested the new algorithm in HMC and MD simulations of unconstrained and constrained dynamics. The tests demonstrated the superiority of the novel scheme over Verlet, BCSS and the HOH-integrator of Predescu et al. \cite{Predescu}. For a wide range of step sizes and MD trajectory lengths,  AIA outperformed other tested integrating schemes  in  accuracy and  sampling efficiency. The analysis of integrated autocorrelation functions and folding evolution demonstrated, for selected sizes of time steps, that AIA possesses up to 5 times better sampling performance than the other tested schemes.

  The idea proposed here may be extended in a natural way to multiple-time-step (MTS) algorithms such as those based on Reversible multiple time scale molecular dynamics \cite{MTS}, the generalised hybrid Monte Carlo method \cite{MTSEscribano}, the Stochastic, resonance-free multiple time-step algorithm \cite{MTSLeimkuhler}, etc. Such extensions are the subject of ongoing work \cite{preprint}.

  In summary, the proposed Adaptive Integration Approach introduces a rational control on integrating the equations of motions in molecular dynamics simulations, leading to enhanced accuracy and performance. To our knowledge this feature was desired but missing by the molecular simulation community.

\section*{Acknowledgements}
  The authors would like to thank the financial support from MTM2013-46553-C3-1-P funded by MINECO (Spain).
  This work has been possible thanks to the support of the computing infrastructure of the i2BASQUE academic network, the technical and human support provided by IZO-SGI SGIker of UPV/EHU and European funding (ERDF and ESF), and the in-house BCAM-MSLMS group's cluster Monako.
  MFP would like to thank the Spanish Ministry of Economy and Competitiveness for funding through the fellowship BES-2014-068640.
  This research is also supported by the Basque Government through the BERC 2014-2017 program and by the Spanish Ministry of Economy and Competitiveness MINECO: BCAM Severo Ochoa accreditation SEV-2013-0323.
  We thank G.A. Papoian (University of Maryland, USA) for valuable discussions.

\appendix\section{Toxin data at $\Delta t / nr = 25$~fs}\label{appendix}
The data collected from the HMC and MD simulations of the unconstrained toxin system at the step size, identified as the stability limit for the velocity Verlet integrator, are presented in \autoref{tab:HMC} and \autoref{tab:MD} respectively.

\noindent
  \begin{table}[h]
  \begin{tabularx}{\textwidth}{XXXXXXX}
    \hline
    HMC Properties & Acceptance Rate (\%) & IACF equilibration & IACF production & Acceptance Rate (\%) & IACF equilibration & IACF production \\
    \hline
    $\Delta t \times L$ (ps) & \multicolumn{3}{|c|}{25 fs $\times$ 1000 = 25 ps}     & \multicolumn{3}{c}{25 fs $\times$ 2000 = 50 ps} \\
    \hline
    VV                       & 42.92               & 105.50         & 102.76         & 39.34               & 226.57         & FAILED   \\
    BCSS                     & 17.64               & 259.38         & FAILED         & 6.01                & 373.43         & FAILED   \\
    AIA                      & 41.30               & 107.04         & 105.94         & 41.10               & 101.02         & 141.67   \\
    \hline
  \end{tabularx}
  \caption{Acceptance rate and IACF for equilibration and production with HMC for the biggest simulated time step and two choices of lengths $L$ of MD trajectories. }\label{tab:HMC}
\end{table}

\noindent
\begin{table}[h]
\begin{tabularx}{\textwidth}{XXX}
    \hline
    MD Properties & IACF equilibration & IACF production \\
    \hline
    $\Delta t/$number of stages (fs) & \multicolumn{2}{|c}{25 fs} \\
    \hline
    VV                               & 124.21         & 122.49    \\
    BCSS                             & 324.56         & FAILED    \\
    AIA                              & 112.75         & 110.98    \\
    \hline
\end{tabularx}
\caption{IACF for equilibration and prodcution with MD for the biggest simulated time step.}\label{tab:MD}
\end{table}

\end{document}